# A Comprehensive Review of the Mapper Algorithm, a Topological Data Analysis Technique, and Its Applications Across Various Fields (2007-2025)

Vine Nwabuisi Madukpe[1,3], Bright Chukwuma Ugoala[2,3], Nur Fariha Syaqina Zulkepli[1]*

[1]School of Mathematical Sciences, Universiti Sains Malaysia, 11800 USM, Penang, Malaysia.
[2]Department of Biostatistics, Augusta University,1120 15th Street, Rains Hall, Augusta, GA 30912, United States.
[3]Department of Mathematics, Abia State University, P.M.B. 2000, Uturu, Abia State, Nigeria.

*Correspondence: farihasyaqina@usm.my Tel.: +6012-543-2006

## Abstract

The Mapper algorithm, a technique within topological data analysis (TDA), constructs a simplified graphical representation of high-dimensional data to uncover its underlying shape and structural patterns. The algorithm has attracted significant attention from researchers and has been applied across various disciplines. However, to the best of the authors' knowledge, no comprehensive review currently exists on the Mapper algorithm and its variants as applied across different fields of study between 2007 and 2025. This review addresses this gap and serves as a valuable resource for researchers and practitioners aiming to apply or advance the algorithm. The reviewed literature comprises peer-reviewed articles retrieved from major academic databases, including Google Scholar, Web of Science, Scopus, JSTOR, PubMed, and IEEE Xplore, using the keywords "topological data analysis", "Mapper algorithm", and "topological graph". The study further provides an overview and a comparative analysis of the suitability of the most commonly used filter functions and clustering algorithms within the Mapper framework. Additionally, it examines current trends, identifies limitations, and proposes future research directions for the Mapper algorithm and its variants, emphasizing the need for developing effective methodologies to streamline the analysis of high-dimensional data in the age of big data proliferation.

## 1 Introduction

The exponential growth of data generated across various disciplines poses significant challenges in analysis due to their complexity and high dimensionality (Faaique, 2023). To address these challenges, innovative methodologies such as topological data analysis (TDA) have become essential for visualizing hidden patterns and extracting valuable information inherent in such datasets. By applying principles of algebraic topology, TDA employs two main techniques, Persistent Homology and the Mapper algorithm, to establish a qualitative framework for analyzing datasets based on their shape, reinforcing the concept that "data has shape, and shape has meaning" (Todd and Petrov, 2022). TDA's growing application across diverse fields has introduced innovative perspectives and approaches to data analysis (Corcoran & Jones, 2023). These TDA techniques enhance the understanding of datasets' inherent structures by extracting topological features and simplifying the representation of high-dimensional datasets (Musa et al., 2021; Chazal & Michel, 2021; Darling et al., 2022).

However, Persistent homology and the Mapper algorithm, as key tools in TDA, serve distinct purposes. Persistent homology is primarily used to extract topological features such as connected components, loops, and voids across multiple scales within high-dimensional data (Shah et al., 2025). It tracks the birth and death of these topological features, presenting them in the form of barcodes or persistence diagrams (Zulkepli et al., 2022). In contrast, the Mapper algorithm is used for visualizing and understanding the underlying topological structure of high-dimensional data. It simplifies high-dimensional datasets into topological graphs by applying a filter function for dimensionality reduction, partitioning the data into overlapping bins, clustering data within each bin, and connecting clusters that share common data points (Carlsson, 2020). This graphical representation provides an effective way to explore complex dataset structures, making Mapper particularly valuable for data exploration and visualization.

The Mapper algorithm has garnered considerable attention in fields such as medicine, biology, finance, social sciences, and other domains where complex, high-dimensional data are prevalent. In biology, it has been employed to analyze gene expression patterns, offering meaningful insights into disease progression, cell differentiation, and protein interactions by uncovering hidden structures within large-scale genomic data (Derwae et al., 2023). In finance, it aids in identifying market regimes, detecting anomalous behaviors, and visualizing risk clusters, thereby providing essential tools for portfolio management and fraud detection (Shiraj et al., 2024). Its applications extend to environmental science (Madukpe et al., 2025), machine learning (Samani & Banerjee, 2023), epidemiology (Ling et al., 2024), biomedicine (Loughrey et al., 2024), and neuroimaging (Singer et al., 2024). In the social sciences, the Mapper algorithm enables researchers to explore and interpret complex social networks, revealing community structures and behavioral patterns that are challenging to detect using traditional statistical methods.

Acknowledging this growing interest, researchers have introduced variants of the Mapper algorithm to address some of its limitations. Among these advancements is the Ball Mapper algorithm, introduced by Dłotko (2019a), which is inspired by the Mapper algorithm. It addresses the use of multiple parameters within the Mapper framework by utilizing only a single parameter in its algorithm. Ball Mapper has found applications across various fields, including financial analysis (Dłotko et al., 2021; Dłotko et al., 2022), biology (Valerio et al., 2024), environmental science (Madukpe et al., 2025), and epidemiology (Dłotko and Rudkin, 2020). In addition to Ball Mapper, several variants of the traditional Mapper algorithm have been proposed, each retaining Mapper principles while introducing specific modifications and enhancements to address its limitations. Notable examples include F-Mapper (Bui et al., 2020), V-Mapper (Imoto and Hiraoka, 2023), G-Mapper (Alvarado et al., 2023), the Mapper-Enhanced Persistent Homology Classification Algorithm (De Lara and Garcia, 2024), Ensemble-based Mapper (Fitzpatrick et al., 2023), and Two-Tier Mapper (Jeitziner et al., 2019). These variants collectively demonstrate the evolution of Mapper techniques, broadening their applicability and enhancing their performance across diverse domains.

Finally, this work aims to provide a comprehensive review of the Mapper algorithm, with a focus on its recent methodological advancements and diverse applications across multiple disciplines. The review is based on peer-reviewed articles published between 2007 and 2025. Through a systematic analysis of these articles, this study aims to demonstrate how the Mapper algorithm has evolved to address its initial limitations by incorporating various algorithmic improvements and introducing emerging variants. The structure of the paper is organized as follows: Section 2.2 presents the list of articles reviewed; Section 2.4 provides an overview of the Mapper algorithm; Section 2.6 offers a comparative analysis of commonly used clustering algorithms and filter functions, emphasizing their effectiveness in preserving dataset structure through Mapper topological graphs; Sections 3.1 and 3.2 present the applications of the Mapper algorithm and its variants in fields such as biology, finance, medicine, environmental science, and machine learning; and Sections 3.3 and 3.4 discuss current challenges, potential solutions, and future directions for improving the scalability, interpretability, and practical utility of the Mapper algorithm and its variants.

# 2 Review Background

Fig. 1 presents a flowchart illustrating the structured approach followed in conducting the review of the Mapper algorithm. It begins with a keyword-based search across multiple academic databases. The collected published articles are categorized into two groups: studies focused on the traditional TDA Mapper algorithm and those discussing its variants. Relevant information is then extracted, such as publication year, author's country, parameters used in the Mapper algorithm, the field of study, data types, and findings, including limitations and future directions. The final stage consists of a discussion and conclusion that reviews the applications of both traditional and variant Mapper algorithms, identifies research challenges, and suggests future directions for further research. This systematic process ensures a comprehensive literature review and evaluation of the Mapper algorithm applications.

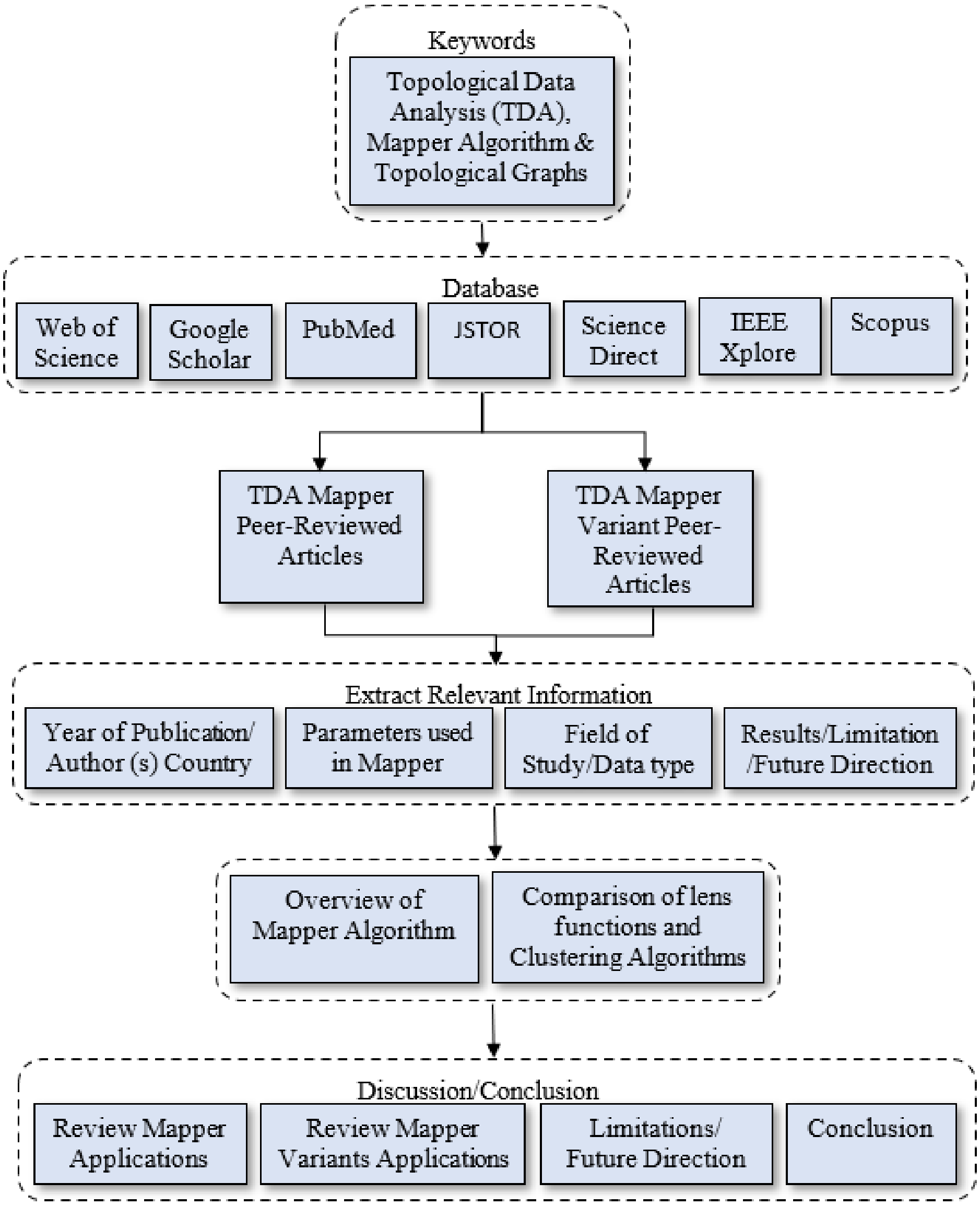

**Fig.1** Methodology Flowchart

## 2.1 List of Abbreviations

| | |
|---|---|
| CDF | Conditional Density Filter |
| FA | Factor Analysis |
| SVD | Singular Value Decomposition |
| L1-∞ | L1-infinity Norm / Regularization |
| UMAP | Uniform Manifold Approximation and Projection |
| DM | Diffusion Maps |
| FCF | Functional connectivity-based features |
| CHL | Chroma and Hue-based Lens |
| Dist-Mean | Mean Distance |
| ISOMAP | Isometric Mapping |
| TAD | Total Absolute Deviation |
| HC | Heteroblasty Classification |
| TDA-GC | Topological Data Analysis Graph Clustering (e.g, greedy modularity) |
| TTMap | Two-Tier Mapper |
| LPA | Label Propagation Algorithm |
| GMM | Gaussian Mixture Models |
| DBSCAN | Density-Based Spatial Clustering of Applications with Noise |
| t-SNE | t-Distributed Stochastic Neighbor Embedding |
| PCA | Principal Component Analysis |
| HACA | Hierarchical Agglomerative Clustering Algorithm |
| MDS | Multidimensional Scaling |
| bht-SNE | Barnes-Hut t-Distributed Stochastic Neighbor Embedding |
| SVD | Singular Value Decomposition |
| SVM | Support Vector Machine |
| KDE | Kernel Density Estimation |
| RLCF | Random Linear Combination of Features |
| MEPHCA | Mapper-Enhanced Persistent Homology Classification Algorithm |
| N/A | Not Available/Not Applicable (when the parameters used are not explicitly stated) |

## 2.2 Articles Screened and Selected for Review

As illustrated in the flowchart, a total of 105 articles, comprising both open-access and non-open-access publications, were initially retrieved from six major academic databases using the keywords "Topological Data Analysis", "Mapper Algorithm", and "Topological Graph". In the preprocessing stage, articles that mentioned these keywords but did not directly apply the Mapper algorithm or its variants in their analysis were excluded. After this filtering process, 89 articles were deemed most relevant. These are organized for detailed review in Table 1a, which covers the Mapper algorithm, and Table 1b, which focuses on its variants.

**Table 1a** List of articles on the Mapper algorithm

| References | Lens/Filter Function | Clustering Algorithm | Data set |
|---|---|---|---|
| Singh et al. (2007) | N/A | N/A | 3D Shape Database |
| Yao et al. (2009) | N/A | N/A | Simulated data on RNA hairpin with GCAA tetraloop |
| Nicolau et al. (2011) | PCA | HACA | high-throughput microarray |
| Lum et al. (2013) | SVD/L1 | HACA | Gene expression/members of the US House of Representatives/ characteristics of NBA basketball players via their performance statistics |
| Pearson (2013) | PCA | HACA | Diabetes patient data with multiple clinical variables. |
| Kyeong et al. (2015) | Correlation | HACA | Resting state fMRI |
| Nielson et al. (2015) | SVD | TDA-GC | preclinical traumatic brain injury (TBI) and spinal cord injury (SCI) |
| Carrière et al. (2017) | N/A | N/A | N/A |
| Guo and Banerjee (2017) | MDS | DBSCAN | chemical process yield prediction and semiconductor wafer fault detection |
| Munch (2017) | N/A | HACA | Simulated |
| Pedoia et al. (2017) | PCA | TDA-GC | MRI data |
| Savic et al. (2017) | N/A | HACA | Physical and chemical properties |
| Mathews et al. (2018) | DM | N/A | RNA-Seq gene expression profiles from The Cancer Genome Atlas (TCGA) |
| Pirashvili et al. (2018) | PCA/MDS/ t-SNE | TDA-GC | 3663 molecules, given in Sybyl line notation |
| Saggar et al. (2018) | t-SNE | HACA | multitask fMRI datasets |
| Siddiqui et al. (2018) | FA | TDA-GC | bronchial biopsy samples and clinical metadata patients with asthma/healthy |
| Robles et al. (2018) | N/A | N/A | 2D grayscale and color images. |
| Campbell et al. (2019) | N/A | HACA | Physiological signals, 85 subjects in response to heat-induced pain. |
| Cawi et al. (2019) | PCA | HACA | Hospital Readmissions and Credit Evaluation |
| Duman et al. (2019) | N/A | HACA | fMRI electroencephalogram (EEG) |
| Geniesse et al. (2019) | PCA | N/A | EEG, fMRI, or MEG data |
| Liao et al. (2019) | | | |
| Patania et al. (2019) | PCA | DBSCAN | N/A |
| Rafique and Mir (2019) | PCA | DBSCAN | Resting-State Functional MRI (rs-fMRI) |
| Wang et al. (2019) | bht-SNE | HACA | Gene expression and survival data |
| Amézquita et al. (2020) | Eccentricity | HACA | scRNA-seq of melanoma tumor cells |
| Belchi et al. (2020) | PCA | DBSCAN | RNA-seq data from lung adenocarcinoma samples |
| Dagliati et al. (2020) | PCA | HACA | Theoritical |
| Qiu et al. (2020) | SVD/L1-∞ | HACA | characteristics from cross-sectional record-level data, not ordered in time |
| Riihimäki et al. (2020) | N/A | N/A | Financial ratios and bankruptcy indicators |
| Tokodi et al. (2020) | PCA | HACA | Branch data obtained from laser scanning of botanical trees. |
| Walsh et al. (2020) | MDS | N/A | |
| Aljanobi et al. (2021) | PCA | HACA | High-dimensional genomic and clinical data. |
| Carr et al. (2021) | N/A | N/A | N/A |
| Endo and Yokogawa (2021) | PCA | DBSCAN | Genome-Based Therapeutic Drugs for Depression |
| Hwang et al. (2021) | t-SNE/ Dist_Mean | DBSCAN | Time-series sensor data |
| Kitanishi et al. (2021) | MDS | N/A | Coronary CTA plaque |
| Lukmanov et al. (2021) | N/A | HACA | 10,188 new depression patients in the JMDC |
| Lopez et al. (2022) | KDE | DBSCAN | Reaven and Miller diabetes dataset |
| Saggar et al. (2022a) | N/A | N/A | Electronic health record (EHR) was conducted on 798 AMI subjects |
| Saggar et al. (2022b) | PCA | HACA | fMRI Brain Image |
| Skaf and Laubenbacher (2022) | MDS | HACA | Two runs of 225 multi-echo EPI images |
| Zhang et al. (2022) | N/A | N/A | Electronic Health Records, genomic sequencing data. |
| Anderson et al. (2022) | | | Functional MRI (fMRI) Dataset |
| Amézquita et al. (2023) | N/A | DBSCAN | Serbia 2014–2015 Multiple Indicator Cluster Surveys (MICS) dataset |
| Chao et al. (2023) | N/A | N/A | RNA-seq data from tumor and healthy subjects. |
| Lauric et al. (2023) | L1-∞ | GMM | Patients who underwent transcatheter edge-to-edge repair (TEER) |
| Derwae et al. (2023) | N/A | N/A | 216 bifurcation aneurysms (90 ruptured) obtained from 3D rotational angiography scans |
| Restrepo (2023) | UMAP | DBSCAN | Mass Spectrometry Imaging (MSI) data from healthy mouse pancreas tissues |
| Duman and Tatar (2023) | N/A | N/A | N/A |
| Escolar et al. (2023) | t-SNE | HACA | Human non-invasive resting state/task Magnetoencephalography (MEG) |
| Fan et al. (2023) | PCA | HACA | 333 major frms' patent portfolios in 1976{2005 across 430 technological areas. |
| Ma et al. (2023) | ISOMAP | | functional MRI (fMRI) time-series data |
| Palande et al. (2023) | N/A | N/A | Parkinson's Progression Markers Initiative (PPMI) |
| Zhou et al. (2023) | PCA/t-SNE | DBSCAN | Vast number of gene expression in plants |
| Rather and Chachoo (2023) | N/A | HACA | CIFAR-10 dataset |
| Ling et al. (2024) | UMAP | HACA | mRNA expression |
| Bungula and Darcy (2024) | SVM | DBSCAN | COVID-19 daily confirmed cases, deaths, and vaccinations |
| Daou and Hanna (2024) | FCF | DBSCAN | Point Cloud |
| Loughrey et al. (2024) | UMAP | LPA | Protein-Protein Interaction (PPI) Datasets |
| Haşegan et al. (2024) | N/A | N/A | Gene expression cancer datasets |
| Hernández-Lemus et al. (2024) | PCA/MDS/ t-SNE | HACA/DBS CAN | Functional magnetic Resonance Imaging (fMRI) datasets |
| Uray et al. (2024) | N/A | N/A | Cardiovascular Signals ECG Data and Heart Rate Signals |
| Percival et al. (2024) | N/A | HACA | Industry 4.0 |
| Samani and Banerjee (2024) | HC | DBSCAN | Leaves of vining grapevine/ maracuya´species |
| Singer et al. (2024) | CHL | N/A | RGB-D images from real and synthetic datasets. |
| Vannoni et al. (2024) | N/A | HACA | Functional magnetic resonance imaging (fMRI) scans |
| Restrepo et al (2024) | N/A | N/A | N/A |
| Zulkepli et al. (2024) | N/A | N/A | N/A |
| Singh et al. (2025) | PCA | HACA | $PM_{10}$ Dataset |
| D. Kim et al. (2025) | N/A | N/A | N/A |
| Madukpe et al. (2025) | PCA | DBSCAN | Air quality data |
| Loughrey et al. (2025) | PCA | HACA | Water quality data. |

**Table 1b** List of articles on the Mapper algorithm variants

| Reference | Mapper Variant Name | Dataset |
|---|---|---|
| Dłotko (2019a) | Ball Mapper | N/A |
| Dłotko (2019b) | Ball Mapper | |
| Dłotko & Rudkin (2020) | Ball Mapper | COVID-19 clinical data |
| Dłotko et al. (2021) | Ball Mapper | N/A |
| Dłotko et al. (2022) | Ball Mapper | N/A |
| Valerio et al. (2024) | Ball Mapper | Micro-CT (micro-computed tomography) images of otolith structures |
| Rudkin et al. (2024) | Ball Mapper | Bitcoin daily return data |
| Madukpe et al. (2025) | Ball Mapper | Air quality data |
| Park et al. (2025) | Ball Mapper | Diffuse large B-cell lymphoma (DLBCL) |
| Bui et al. (2019) | F-Mapper | NKI breast cancer, |
| Jeitziner et al. (2019) | Two-Tier Mapper | Microarray and RNA-seq data |
| Kang and Lim (2021) | Ensemble-based Mapper | Diabetic data |
| Fitzpatrick et al. (2023) | Ensemble-based Mapper | TCGA-BRCA (breast cancer) gene expression |
| De Lara and Garcia (2024) | MEPHCA | N/A |
| Imoto and Hiraoka (2023) | V-Mapper | single-cell gene expression data of endocrine cells |
| Alvarado et al. (2023) | G-Mapper | N/A |
| Tao and Ge (2025) | D-Mapper | Synthetic datasets and SARS-CoV-2 RNA seq |

## 2.3 TDA Mapper Research Trend

The Mapper algorithm has gained significant global recognition, as illustrated in Fig. 2**.** Fig. 2a depicts the increasing number of publications utilizing the algorithm over time. From 2007 to 2015, research activity remained minimal, with only a slight increase. However, after 2017, interest in the algorithm grew significantly, culminating in a sharp rise in 2019. This upward trend continued, peaking in 2024, reflecting the growing adoption of the Mapper algorithm. Data for 2025 is still being collected and may evolve as the year progresses. Overall, this trend highlights the algorithm's increasing significance in academic and applied research. Meanwhile, Fig. 2b presents the geographical distribution of authors contributing to Mapper-related studies. The USA leads by a substantial margin, highlighting its dominant role in the algorithm's development and application. The UK follows with strong contributions, while Japan, China, Switzerland, France, and Malaysia exhibit moderate research activity in this area. Although contributions from other countries indicate a growing international presence, the Mapper research remains concentrated in a few leading research hubs. While global diversification is expanding, it remains relatively limited.

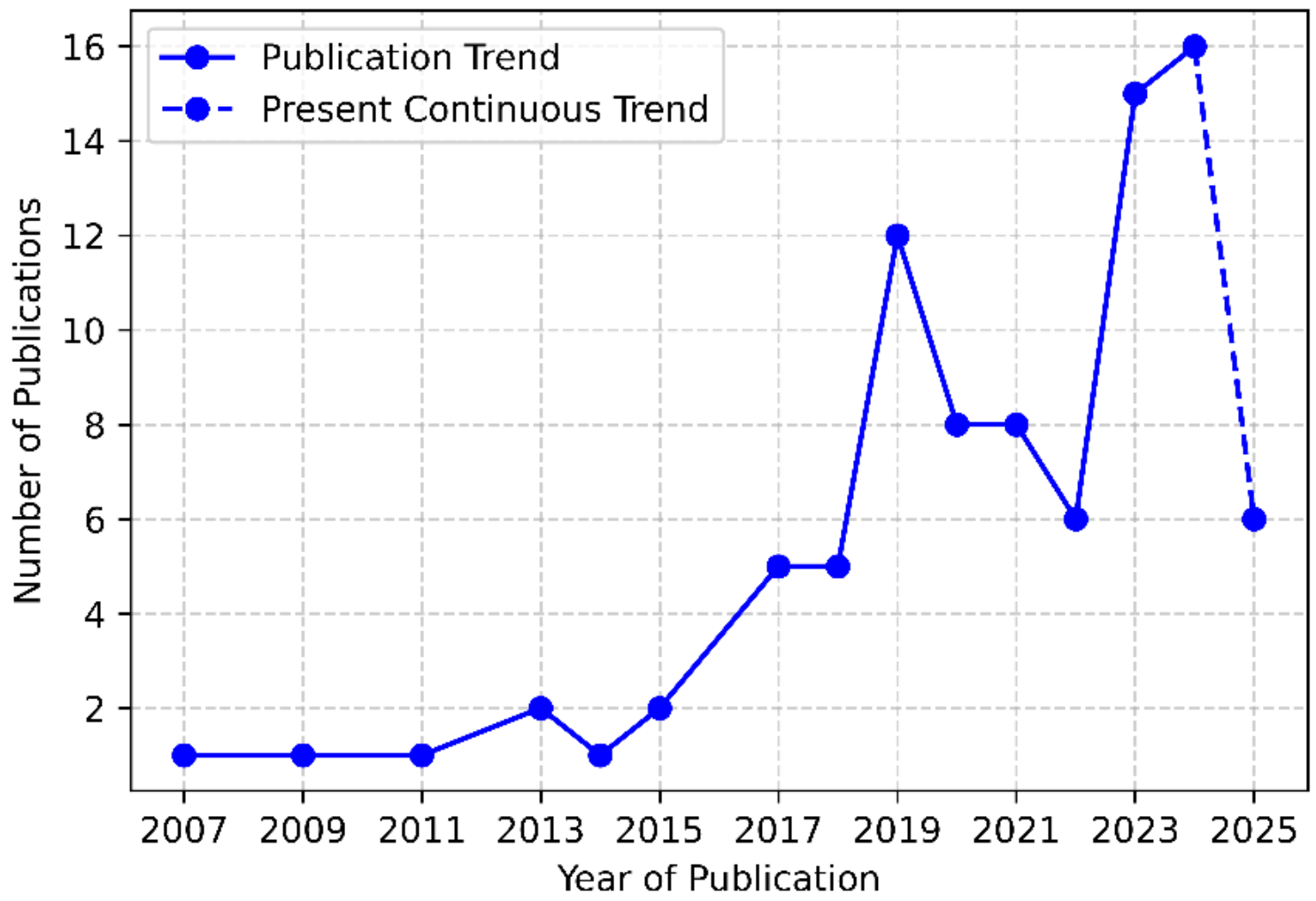

(a)

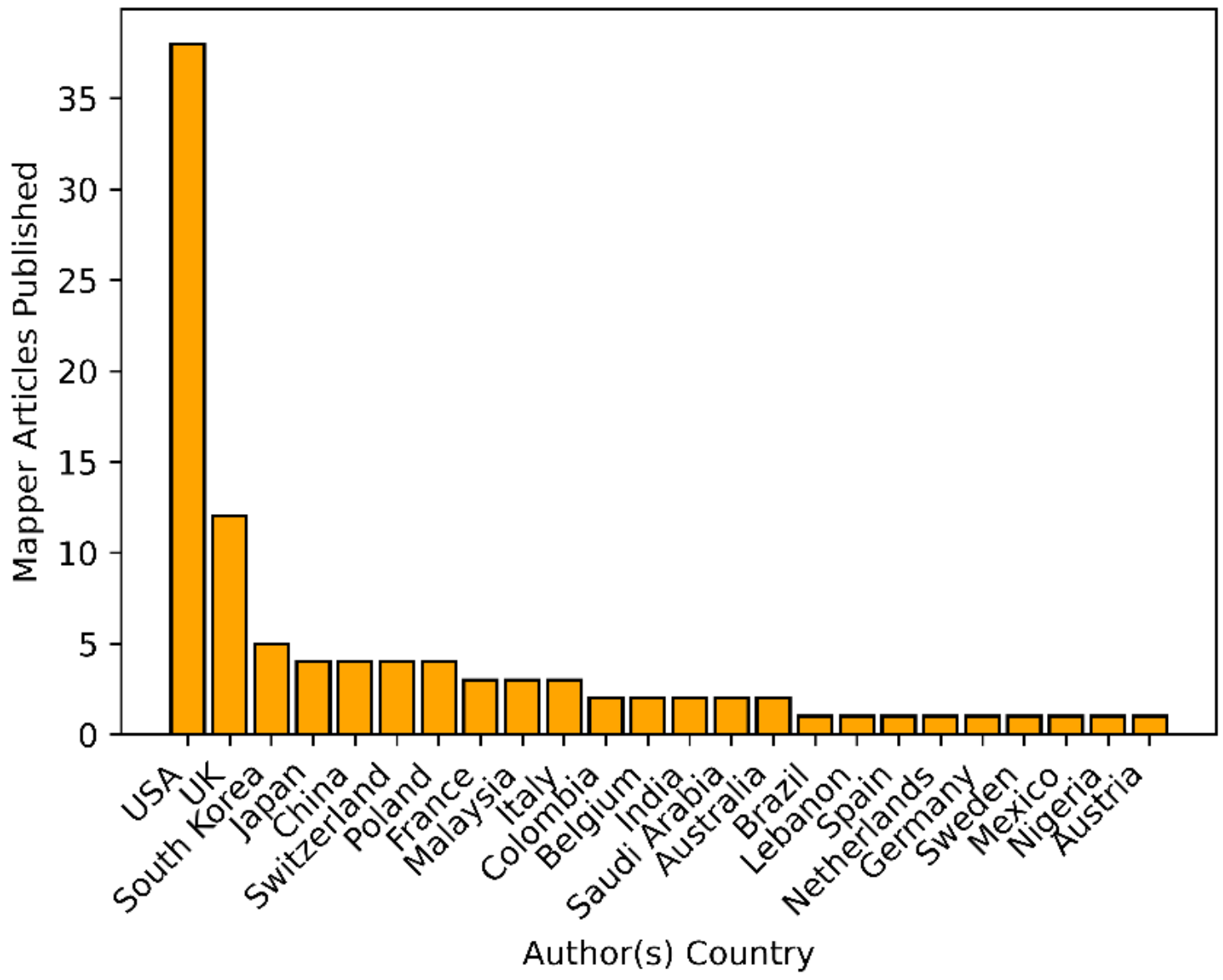

(b)

**Fig. 2** Mapper Algorithm trends for (a) publication and (b) geographical distribution of authors

## 2.4 Overview of Mapper Algorithm

The Mapper algorithm was first introduced by Singh et al. (2007) and further developed by Nicolau et al. (2011). Over time, it has gained significant attention and has been successfully applied across various domains, as highlighted by Lum et al. (2013). It is a pivotal method within TDA. It offers an approach to visualizing the topological structure of high-dimensional datasets or point clouds by employing multiple parameters, such as a filter function, cover, overlap percentage, and clustering algorithm, to produce a simpler dataset representation. The entire algorithm is fundamentally built on partially clustering the data, guided by a collection of functions defined over it. The following sections outline the fundamentals on which the Mapper algorithm is built.

### *2.4.1 Reeb Graph*

Given a topological space $M$ and a continuous real-valued function $f : X \to M$, the Reeb graph identifies the level set of $f$ at value $c$, $f^{-1}(c) = \{m \in M \mid f(m) = c\}$, representing all points with the same function values. Each connected component of a level set is collapsed into single points forming nodes. The edges connecting the nodes represent how the components merge or split as the function value changes (Carrière et al., 2017). Fig. 3 is the construction of a Reeb graph from a Torus. When the target parameter space is □ , the construction represents the Mapper graph, a stochastic version of the Reeb graph because the Mapper graph approximates the Reeb graph by clustering data points based on a filter function and then connecting these clusters as discussed in Section 2.3.2. If the covering □ is too broad, the result approximates the Reeb graph's image. Conversely, with a sufficiently refined cover and appropriate clustering, the Mapper construction can closely approximate, or even accurately recover, the Reeb graph structure, as demonstrated by Singh et al. (2007) and further formalized by Carrière et al. (2017).

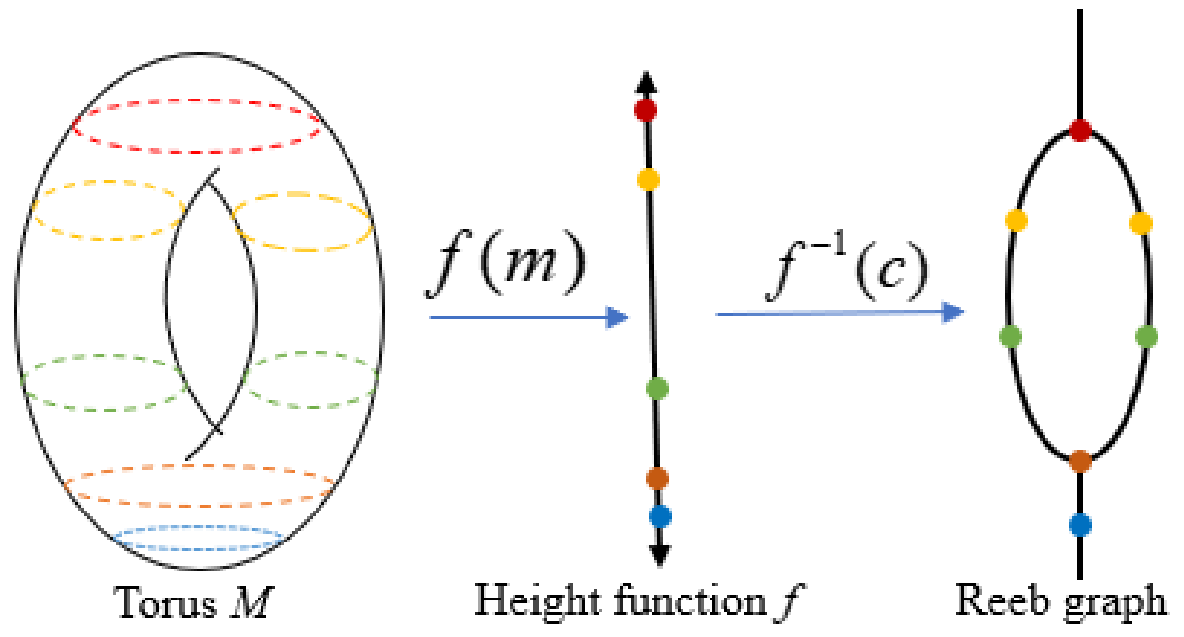

**Fig. 3** Reeb graph construction for a Torus using a height function

### 2.4.2 *The Topological Version of Mapper*

Let $X$ and $X^*$ be topological spaces and $f: X \rightarrow X^*$ be a continuous map. Let $U = \{U_i\}_{i\in I}$ be a finite open set covering $X^*$ then consider $f^*(U) = \{f^{-1}(U_i)\}_{i\in I}$ the covering of $X$ where $\{f^{-1}(U_i)\}_{i\in I}$ is called the pullback of $U$ along $f$. The Mapper topological graph $G$ Construction is defined to be the nerve complex of the pullbacks covering (Singh et al. 2007).

$$G(U,f) = N(f^*(U)) = \left\{\sigma \subseteq f^*(U) : \bigcap_{f^{-1}(U_i)\in\sigma} f^{-1}(U_i) \neq \varnothing\right\} \quad (1)$$

$\sigma$ represents a subset of the pullback cover $f^*(U_i)$

A simplex $\sigma$ is included in the nerve (collection of sets) $N(f^*(U))$, if the interception of all the pre-image $f^{-1}(U_i)$ in $\sigma$ is non-empty, forming a connected component of nodes and edges, if otherwise singleton nodes are formed. This nerve complex captures the interconnection of the pullback cover, representing a relationship between subsets of $X$ induced by the continuous map $f$ and the original cover $U$ on $X^*$ (Bui et al. 2020).

### 2.4.3 *The Statistical Version of Mapper*

In constructing a Mapper topological graph, four steps are followed with parameter selection and setup, as shown in Fig. 4. The parameters, such as the filter function, the clustering algorithm, the setup of the cover, and the percentage overlap of the covers, are determined (Singh et al., 2007). The filter function is used to discriminate the characteristics of the data. Next, the filtered range is organized into overlapping bins called covers. A specified overlap percentage determines how much these covers overlap. A clustering algorithm is then applied to the pullback of the data points within each cover based on their proximity. This process produces a topological graph, a simplified representation of a high-dimensional dataset. In this graph, nodes represent the clustered data points, and edges indicate nodes that share common data points based on the specified overlap percentage.

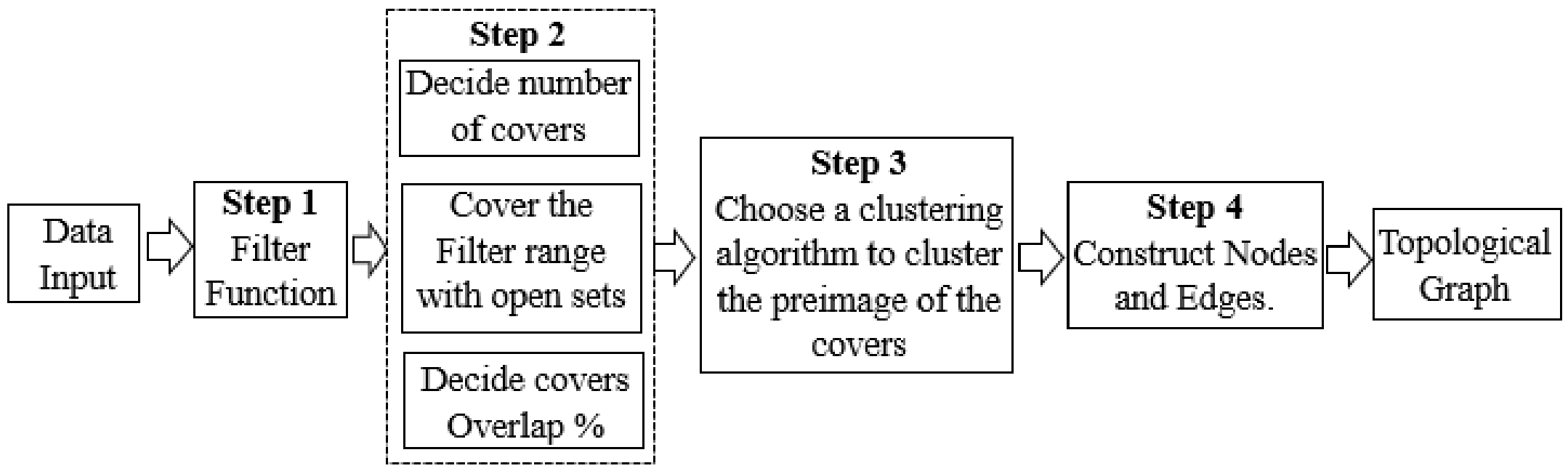

**Fig. 4** Mapper Flowchart

Given a data set or a point cloud $X \subset \mathbb{R}^n$ consisting of a collection of points, the steps for constructing a Mapper topological graph are as follows:

*Step 1*: Choose a filter or lens function, $f : X \to \mathbb{R}^d$ to discriminate the characteristics of the data set by projecting $\mathbb{R}^d$ to a lower Euclidean space, $\mathbb{R}^d$ where $(d < n)$.

Dimension reduction methods such as PCA, t-SNE, UMAP, or custom dimension reduction functions can be used as filter functions. Fig. 5 illustrates the frequency of using various lens/filter functions in the Mapper algorithm across the reviewed articles. PCA is the most commonly used lens function, significantly outpacing other techniques. t-SNE, MDS, and UMAP are also relatively popular, though their usage is much lower than PCA. A long tail of alternative methods, including SVD, L1-Infinity, and Eccentricity, among other custom filter functions, indicates a diverse range of approaches. This distribution suggests that while PCA remains the dominant choice, researchers explore multiple dimensionality reduction and feature extraction methods depending on their dataset and analysis needs.

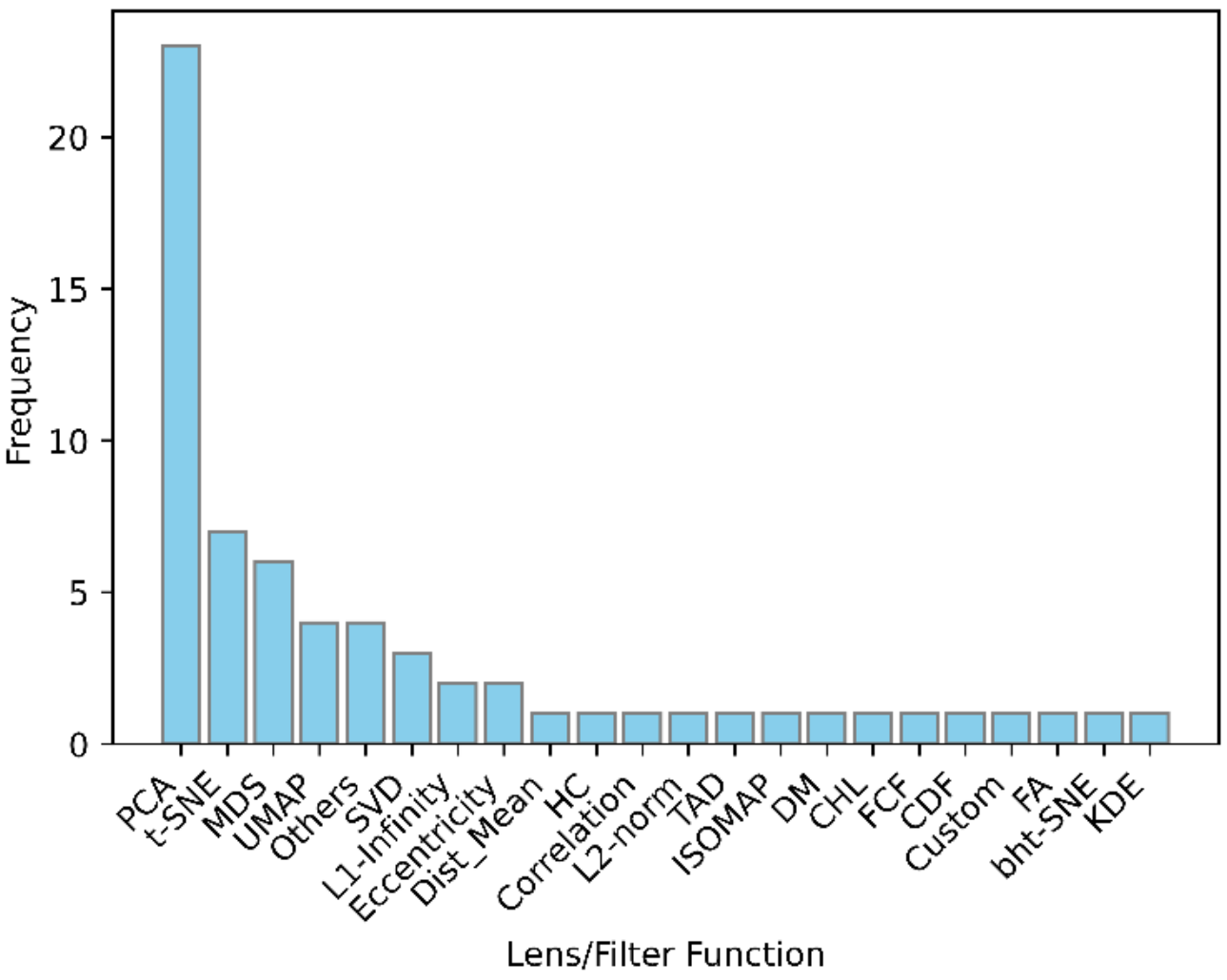

**Fig. 5** Distribution of lens or filter Functions used in Mapper analyses across the reviewed articles

*Step 2*: Build a cover on the image of the discriminated values through a dimension-reduced function $f(X)$, with a collection of open sets, $(U_i)_{i \in I}$ where $I$ is a finite indexing set. Typically, the set $(U_i)_{i \in I}$ is composed of overlapping cubes (products of intervals), and at this stage, it is essential to determine the percentage of overlap.

*Step 3*: For each interval or cover $U_i$, cluster the points in the preimage $f^{-1}(U_i)$ by utilizing a suitable clustering algorithm such as HACA, DBSCAN, or custom clustering algorithm to form a set of clusters $C_{i_n}$.

Fig. 6 presents the frequency of different clustering algorithms used in the Mapper algorithm across the reviewed articles. HACA is the most frequently used, indicating a strong preference for its hierarchical structure and ability to capture data relationships effectively. DBSCAN is the second most common choice, suggesting its utility in identifying clusters of varying densities. TDA-GC is used considerably less often, while other methods such as TTMAP, GMM, LPA, and others are more of a customized approach. This distribution highlights that while HACA and DBSCAN dominate, researchers occasionally experiment with alternative clustering strategies to optimize Mapper analyses.

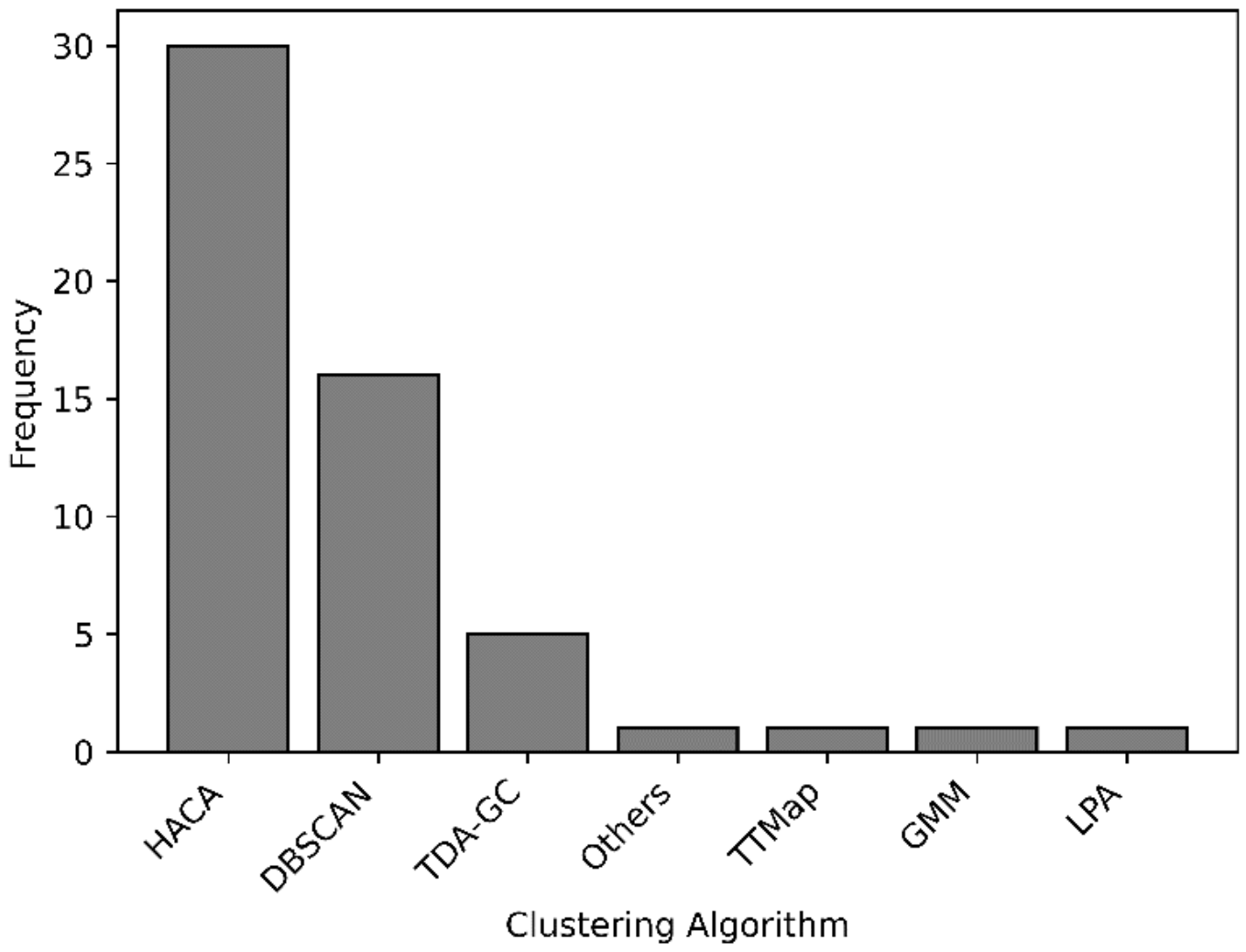

**Fig. 6** Distribution of clustering algorithm used in Mapper analyses across the reviewed articles

*Step 4*: Formation of Mapper topological graph. Each cluster $C_{i_n}$ forms node $(N)$ and for any $i \le j \le I$, if clusters $C_{i_n}$ and $C_{j_p}$ share the same data points and overlap, that is, if $C_{i_n} \cap C_{j_p} \ne \varnothing$ then an edge $(E)$ (a line) is drawn connecting the two nodes; if otherwise, singleton nodes are formed, yielding a topological graph $G = (N, E)$.

## 2.5 Mapper Algorithm Implementations Software

Several software implementations of the Mapper algorithm have been developed to support its application across diverse fields, each offering distinct features tailored to specific research and analytical needs. These implementations provide researchers and data scientists with powerful tools for identifying patterns, clustering complex datasets, and intuitively visualizing high-dimensional structures in the form of a Mapper topological graph. The choice of a particular Mapper implementation depends on the requirements of a given project, whether prioritizing flexibility, ease of use, computational efficiency, or commercial support. The available software implementations and their key features are presented in Table 2, offering a comparative overview to assist users in selecting the most suitable tool for their specific applications.

**Table 2** The software for the Mapper algorithm and Variants is implemented in Python or R

| Software | Advantages | Disadvantages |
|---|---|---|
| KeplerMapper Van Veen et al. (2019) | Open-source and actively maintained. User-friendly with comprehensive documentation. Integrates well with Python data science libraries. | Scalability Issues Requires familiarity with Python. Subjectivity in Visualization |
| Python Mapper Müllner (2013) | Integration with Python Ecosystem Flexible and Customizable | Not Always Scalable for Big Data Limited Documentation and Community Support |
| tda-mapper Simi, L. (2024). | Optimized for large-scale data. | Relatively newer, smaller user base. |
| Mapper Interactive Zhou et al. (2020) | Interactive visualization capabilities. Scalable and extendable for large datasets. | Web-based interfaces may require specific setups. Web-based interfaces may require specific setups. |
| TDA Mapper Pearson (2013) | Integrates seamlessly with R | It may not be as actively maintained. |
| Sakmapper Szairis (2016) | Automated Selection of Parameters It handles large and complex datasets Faster computation | Potential Loss of Detail Less Control for Experts Limited Adoption and Documentation |
| Ayasdi<br>Lum et al. (2013) | User-friendly with robust enterprise support.<br>Advanced analytics features. | Performance limitations with large datasets.<br>Proprietary software requires a license. |
| pyBallMapper<br>Gurnari (2022),<br>Dłotko (2019b) | Open-source and actively maintained.<br>User-friendly with comprehensive documentation.<br>Integrates well with Python data science libraries. | Performance challenges with large datasets.<br>Requires familiarity with Python and R |
| R Package<br>Ballmapper<br>Dlotko (2019c) | Can handle various types of datasets (numerical, mixed types with some preprocessing)<br>Flexible integration in R | Basic visualization in R<br>Smaller user community and less documentation |

Fig. 7 illustrates the implementation process of the Mapper algorithm on a b-shaped point cloud, $X \subset \mathbb{R}^2$ demonstrating how the algorithm generates a Mapper topological graph. Firstly, the point cloud $X$ is discriminated and projected onto a one-dimensional space along the $y-$axis using a filter function $f: X \to X^*$ where $X^* \subset \mathbb{R}$ and subsequently covered with open sets $U_i = \{U_1, U_2, U_3, ..., U_8\}$ with each cover overlapping in an equal percentage, as shown in Steps 1-2. In Step 3, clustering is then performed on the preimages or pullback of $U_i$, $f^{-1}(U_i)$ to have □ , where $U_1$, $U_2$,... corresponds to $C_1, C_2, \ldots$ and if $C_1 \cap C_2 \neq \varnothing$ then an edge $E$ is formed connecting the nodes $N$ formed by $C_1$ and $C_2$. However, if $C_1 \cap C_2 = \varnothing$ then $C_1$ and $C_2$ forms two singleton nodes. The collection of nodes $N$ and edges $E$ In Step 4, we produced a Mapper topological graph $G(N, E)$.

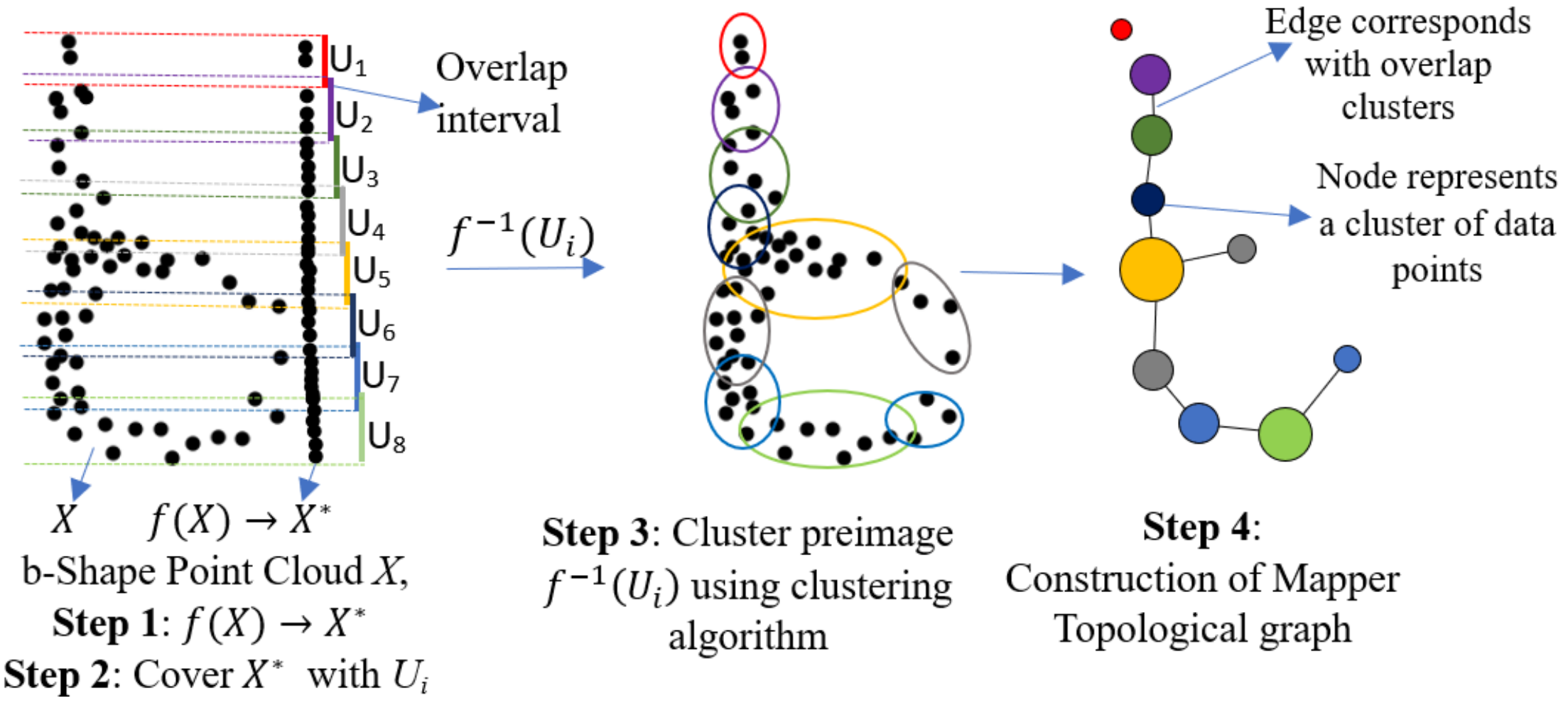

**Fig. 7** Toy b-shape point cloud used in illustrating Mapper algorithm implementation

## 2.6 Comparison of Clustering Algorithms and Lens Functions in Preserving Dataset Topology Using the Mapper Algorithm

Five clustering methods, GMM, TDA-GC, HACA, DBSCAN, and K-Means, were implemented within the Mapper algorithm framework using KeplerMapper. The algorithm was applied to a two concentric circle-shaped point cloud dataset, with the following settings: a projection onto $y-$axis as the filter function, 15 cover partitions, a 20% overlap between covers, and a fixed sample or cluster number of 10 across all clustering methods. The resulting topological graphs are presented in Fig. 8. Among these, GMM and TDA-GC in Fig.8b and 8e best preserved the circular topology of the dataset by maintaining a well-connected structure that accurately reflects the underlying loops. However, they failed to distinguish the two separate rings, as their focus is on global connectivity rather than strict cluster separation. DBSCAN in Fig. 8d, on the other hand, successfully recognized the two distinct circles, making it the most effective method in preserving both topology and structural separation. HACA in Fig. 8c, though appearing fragmented, demonstrated a different strength by capturing fine-grained structures and outliers, making it more sensitive to local variations and anomalies.

While HACA may seem overly disconnected, its ability to identify subtle variations and outliers makes it useful for detailed structural analysis. In contrast, GMM and TDA-GC focus on overall smoothness and connectivity, which is beneficial when studying global geometric patterns. DBSCAN balances between topology and separation, making it ideal for datasets with distinct but continuous structures. K-Means in Fig. 8f, however, performed less, as its assumption of spherical clusters led to a highly scattered and disconnected representation. Ultimately, the best clustering method depends on the analytical goal, whether the focus is on preserving global topology, detecting distinct clusters, or analyzing finer details within the data.

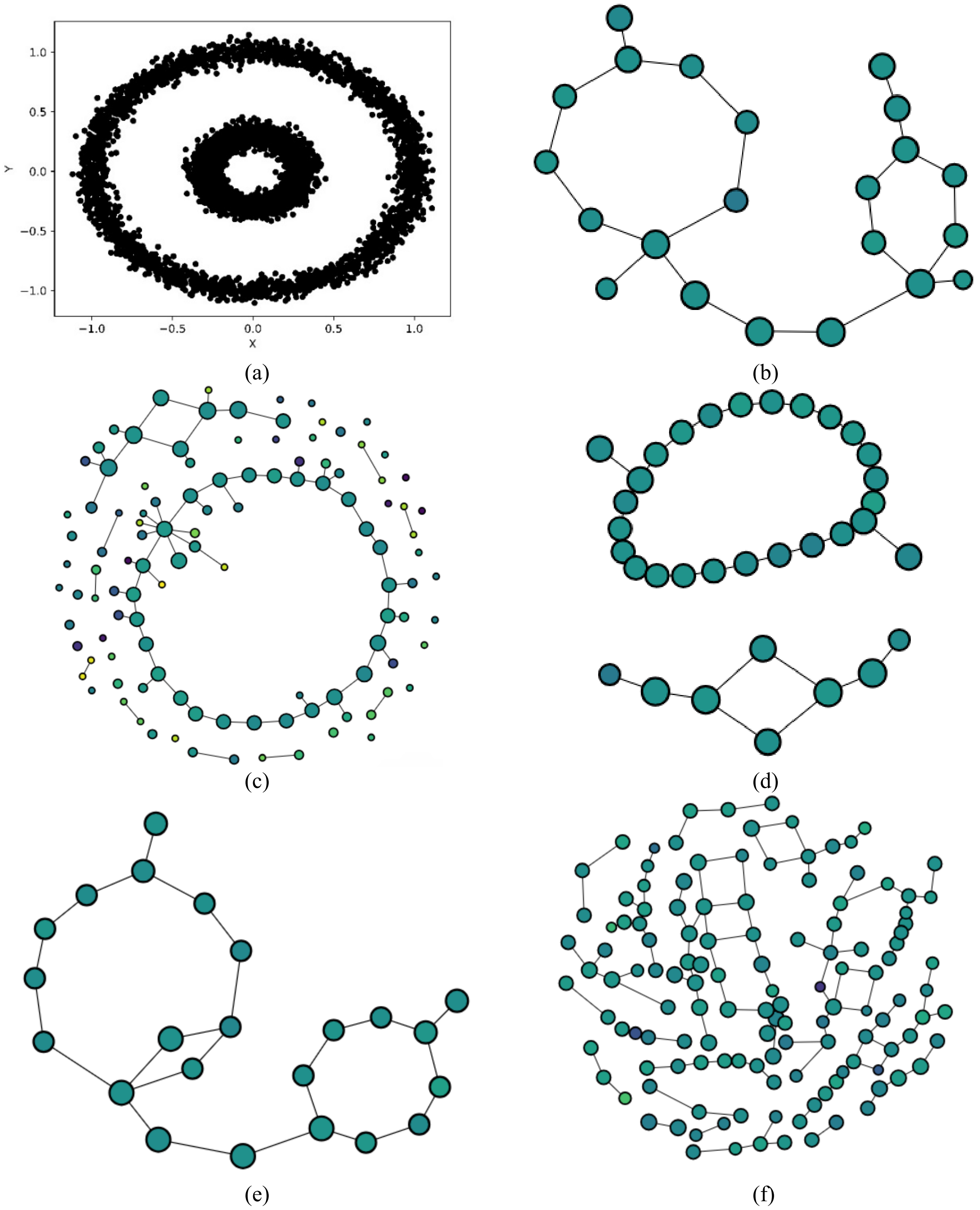

**Fig. 8** Comparison of Clustering Methods in Mapper Algorithm for Optimal Dataset Representation (a) two concentric circle-shaped point cloud (b) GMM (c) HACA (d) DBSCAN (e) TDA-GC (f) k-Mean

Fig. 9 illustrates the evaluation of the performance of different lens functions in the Mapper algorithm for representing the two concentric circle-shaped point cloud datasets in Fig. 8a above. The TDA-GC clustering method is applied using KeplerMapper with 15 cover partitions, 20% overlap, and a fixed sample/cluster number of 10. The four Mapper topological graphs in Fig. 9(a, b, c, and d) are produced by PCA, t-SNE, UMAP, and MDS, respectively. PCA and MDS demonstrate the best performance in preserving global and local structures of the dataset. PCA Fig. 9a effectively retains the global circular structure, accurately representing inner and outer loops with minimal distortion. MDS Fig. 9d also preserves the dataset's global topology while maintaining local connectivity, making it a strong alternative. In contrast, t-SNE Fig. 9d disrupts the global structure, fragmenting clusters and distorting the connectivity between points. UMAP Fig. 9c further exacerbates this issue, breaking the dataset into disconnected components and losing essential local relationships. Therefore, PCA emerges as the most robust lens function, offering the best balance between global topology preservation and local connectivity, with MDS following closely behind, however, this assertion remains hypothetical as it may not reflect the situation for a higher dimensional dataset.

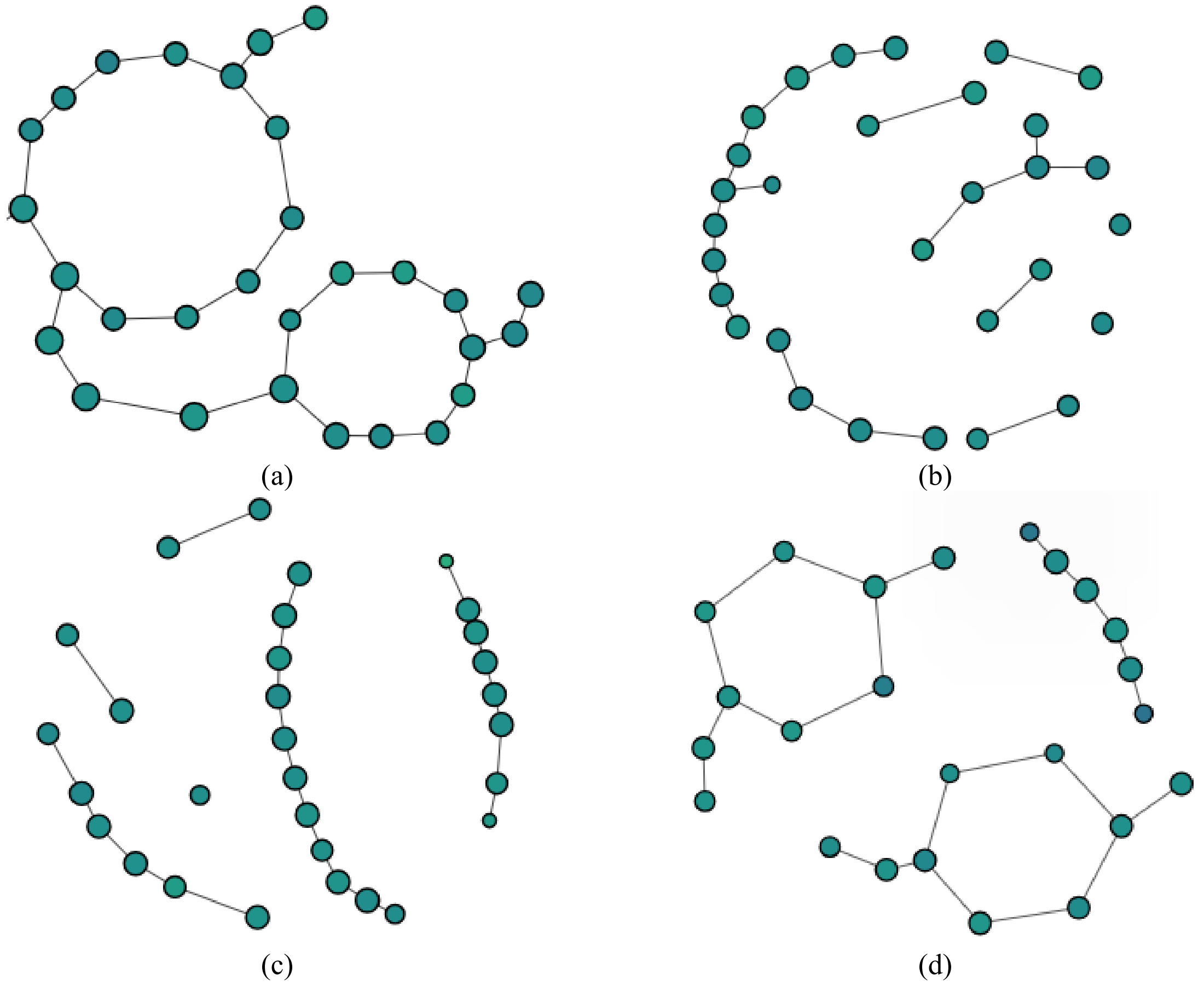

**Fig. 10** Comparison of Lens Function in Mapper Algorithm for Optimal Two Concentric Circle-Shaped Point Cloud Dataset Representation (a) PCA (b) t-SNE (c) UMAP (d) MDS

# 3 Discussion

This section provides a detailed discussion of the Mapper algorithm's evolution, emphasizing its practical applications, algorithmic advancements, and persistent limitations. It describes how Mapper has been implemented across diverse domains, such as neuroscience, environmental science, and machine learning, to derive meaningful insights from high-dimensional data through topological graphs. It further introduces various algorithmic variants developed to enhance Mapper's robustness, scalability, and interpretability. These innovations include methods for reducing sensitivity to parameter choices, improving temporal analysis, and automating cover construction. Moreover, it addresses Mapper's limitations, including its reliance on user-defined parameters and the potential loss of global structural information, while suggesting future research directions focused on refining filter function selection and strengthening theoretical frameworks for more consistent topological insights.

## 3.1 Applications of the Mapper Algorithm

### *3.1.1 Computational Mathematics and Algorithm*

Singh et al. (2007) employed the Mapper algorithm as a computational approach to simplify and visualize high-dimensional datasets while maintaining their topological structure. This work provided a novel way to qualitatively analyze massive datasets, with applications spanning computer graphics and medical research. As interest in Mapper grew, Carrière and Oudot (2017) investigated its theoretical properties, focusing on the stability of 1-dimensional Mapper and its relation to the Reeb graph. Their study provided mathematical guarantees on Mapper's reliability in capturing topological features while minimizing sensitivity to parameter selection. This was further extended by Belchí et al. (2019), who quantified the instability of Mapper-type algorithms, introducing a numerical measure to identify optimal parameter ranges and ensure more consistent outputs. With Mapper gaining traction across scientific disciplines, researchers focused on refining its efficiency and usability.

Munch (2017) provided a foundational introduction to the field, making the field more accessible by formalizing two key TDA methods, such as persistence diagrams for capturing loops and holes in data and Mapper for visualizing data structures through topological graphs. Robles et al. (2018) applied the Mapper algorithm to image analysis, highlighting the limitations of conventional feature extraction techniques that often fail to capture global image structure. Their study demonstrated that Mapper could effectively identify meaningful topological features, offering an advanced approach to shape recognition and classification in image datasets. Jeitziner et al. (2019) contributed to clustering advancements by developing Two-Tier Mapper, an unbiased topology-based clustering method for global gene expression analysis. Traditional clustering methods lacked stability and relied heavily on user-defined parameters, but TTMap integrated hierarchical and partitioning clustering techniques, leading to more robust subgroup identification in biological datasets. Van Veen et al. (2019) introduced Kepler Mapper, a flexible Python implementation of Mapper, addressing the need for an open-source, user-friendly tool that integrates seamlessly with machine learning pipelines, expanding Mapper's accessibility across disciplines.

Walsh et al. (2020) introduced TDAview, an interactive online visualization tool for Mapper-based topological analysis. Their tool addressed the lack of user-friendly TDA software, enabling biologists and clinicians to generate and explore Mapper topological graphs without programming knowledge. More recent studies have focused on refining the Mapper for computational efficiency and real-world applications. Finally, Bungula and Darcy (2024) further examined Mapper's stability in point cloud data, showing that clustering parameter variations, especially in DBSCAN, could affect Mapper's robustness. Their findings demonstrated that increasing cover size alongside clustering parameter adjustments improved Mapper stability, ensuring more reliable topological analyses. De Lara and Garcia (2024) introduced the Mapper-Enhanced Persistent Homology Classification Algorithm (MEPHCA) to refine classification tasks in machine learning. While standard Mapper approaches provided powerful visualizations, they often incurred high computational costs. Their study demonstrated that MEPHCA significantly reduced these costs while maintaining accuracy, advancing the computational robustness of Mapper-based methods.

Vannoni et al. (2024) proposed a method for selecting Mapper algorithm parameters by combining quantitative metrics like Netsimile and graph topological features with qualitative visual analysis. Using brain features from the Human Connectome Project, they identify parameter configurations that improve graph stability and data representation. Their findings emphasize the role of quantitative assessment in parameter selection and support broader adoption through an open-source approach. Restrepo and Cornelis (2024) integrate covering-based rough sets with the Mapper algorithm from TDA to analyze numerical datasets. By adjusting the Mapper's parameters, they generate coverings that define lower approximations and assess classification quality. Their approach enhances rough set analysis, with future work focusing on parameter optimization and broader applications. Collectively, these advancements highlight how Mapper has evolved from a theoretical tool into a powerful computational framework for high-dimensional data analysis. It now offers improved clustering, stability, and visualization tools for diverse scientific and industrial applications, bridging the gap between theory and real-world usability. Loughrey et al. (2025) introduced a hotspot detection algorithm built on the Mapper framework to identify clinically meaningful subgroups in high-dimensional biomedical data. Applied to ER+ breast cancer cohorts, the method revealed a novel patient subgroup with poor prognosis not captured by existing classifications, which was validated in an independent dataset with consistent survival and molecular signatures.

### *3.1.2 Environmental Science*

TDA Mapper has been increasingly applied in ecological and environmental studies, offering new insights where traditional statistical methods fall short. Savic et al. (2017) pioneered the use of TDA in soil science, uncovering intricate relationships between soil properties, land use, and climate. Their study provided a deeper understanding of anthropogenic effects on soil formation, overcoming the limitations of conventional analytical techniques. In a more recent application, Endo & Yokogawa (2022) utilized TDA to analyze indoor carbon dioxide ($CO_2$) concentrations and plug-level electricity usage in an active learning space. By employing the Mapper algorithm, they successfully visualized high-dimensional sensor data, identifying strong correlations between $CO_2$ levels, energy consumption, and human activity patterns, insights that traditional statistical approaches struggled to reveal.

Building on these advancements, Zulkepli et al. (2024) explored spatial patterns of particulate matter ($PM_{10}$) in Malaysia using TDA Mapper. Traditional methods had difficulty visualizing relationships between monitoring stations, but the Mapper algorithm facilitated the construction of a topological graph revealing homogeneous $PM_{10}$ regions. This can help in improving air pollution monitoring and policymaking. Most recently, Madukpe et al. (2025) conducted a comparative study of Ball Mapper and conventional Mapper to assess air pollutant behavior across 60 air quality monitoring stations in Malaysia. Their findings suggested that while both methods provided valuable insights, Ball Mapper's reliance on a single parameter made it more effective for analyzing pollutant trends and detecting redundancies in monitoring networks. These studies highlight the evolving role of TDA in environmental science, demonstrating its ability to address complex, high-dimensional data challenges and improve decision-making processes in ecological and pollution monitoring.

### *3.1.3 Neurosciences*

TDA Mapper has played a transformative role in neuroscience, offering innovative methods to analyze brain dynamics and gene expression networks. Early applications of TDA in this field date back to Kyeong et al. (2015), who used the Healthy State Model (HSM) and Mapper algorithm to analyze functional connectivity data in individuals with Attention-Deficit or Hyperactivity Disorder (ADHD). Their study provided an objective method for assessing ADHD severity by identifying significant correlations between symptom severity and disease-specific neural components through topological graphs. Saggar et al. (2018) further expanded TDA Mapper applications by tracking whole-brain activity maps during multitask environments. Traditional functional connectivity methods struggled to retain spatial and temporal details, but Mapper through topological graphs enabled a fine-grained analysis of transient neural states, successfully predicting cognitive task performance.

In 2019, Patania et al. applied TDA to gene co-expression networks, revealing that Mapper-generated topological graph structures aligned with known anatomical organizations, particularly in dopamine-related pathways. Around the same time, Geniesse et al. (2019) introduced Dynamic Network Structure Reconstruction (DyNeuSR), which uses tools like the Mapper algorithm and persistent homology to extract multi-scale structural features from high-dimensional data. It is an open-source platform designed to extract neurophysiologically valid representations from high-dimensional neuroimaging data without collapsing spatial or temporal information, enabling more precise single-participant analyses. That same year, Duman et al. (2019) investigated dynamic brain reconfiguration during a working memory n-back task using Mapper on magnetoencephalography (MEG) data. Their findings through Mapper topological graph revealed distinct brain state transitions, offering a novel neuroimaging marker for working memory performance. As interest in TDA-based neuroimaging grew, Saggar et al. (2022a) studied how methylphenidate (MPH) influences brain function under cognitive load and anxiety. Their research demonstrated that MPH enhances neural resource engagement under such conditions, improving our understanding of cognitive enhancement therapies.

Building on these advancements, Geniesse et al. (2022) and Saggar et al. (2022b) developed NeuMapper to address computational inefficiencies in traditional Mapper methods, improving parameter selection and enabling better characterization of whole-brain dynamics during multitask cognition. Also, Zhang et al. (2022) introduced the Temporal Mapper, bridging data-driven and mechanistic models of brain function by extracting attractor transition networks from neural time series data. Their work successfully predicted behavioral performance and characterized neural transitions in real time. In 2023, Fan et al. (2023) and Bhad et al. (2023) applied Mapper-generated topological graphs to examine the effects of low-frequency repetitive transcranial magnetic stimulation (rTMS) on brain dynamics. Their studies showed that rTMS significantly modified temporal brain configurations, with distinct effects observed in the frontal and occipital lobes. Lauric et al. (2023) used the Mapper algorithm to classify aneurysm rupture status in 216 bifurcation aneurysms, analyzing 12 shapes and 18 radiomics features. Ruptured aneurysms had higher Mapper dissimilarity scores (MDS), indicating greater shape variability. An MDS threshold of 0.0417 achieved an AUC of 0.73, and incorporating MDS into multivariate analysis improved classification (AUC = 0.82), outperforming shape and radiomics models. Findings support the Mapper potential for aneurysm assessment.

The innovations have enhanced TDA-based methodologies in neuroscience, making them more robust and insightful. Duman and Tatar (2023) leveraged the Mapper algorithm to visualize dynamic fluctuations in neural activity, identifying biomarkers that differentiate cognitive tasks. Haşegan et al. (2024) systematically examined Mapper parameter choices, showing that optimized selection could enhance the extraction of meaningful brain activity transitions from fMRI data. Most recently, Singer et al. (2024) applied TDA Mapper neuroimaging analysis to investigate the effects of psilocybin on meditation. Their findings suggested that psilocybin enhances meditation-induced brain states by altering topological brain structures, offering new insights into potential therapeutic applications. Collectively, these studies highlight the power of TDA in refining neuroimaging analyses, improving computational frameworks, and revealing novel aspects of brain function and cognition.

#### *3.1.4 Machine Learning and Algorithm Optimization*

Pearson (2013) pioneered the use of the Mapper algorithm to enhance the interpretability of artificial neural networks (ANNs), addressing the challenge of visualizing high-dimensional clustering. Traditional clustering methods often struggled to effectively interpret neural decision boundaries, but Pearson demonstrated that Mapper-generated visualizations or topological graphs provided meaningful insights into classification performance. The study also showed that Mapper could track decision boundary evolution in neural networks, improving model validation and refinement. Around the same time, research on ensemble learning for Mapper parameter optimization highlighted the instability caused by varying parameters. The findings suggested that ensemble learning could enhance parameter selection, leading to more robust topological structures in data analysis.

As the field evolved, researchers sought to integrate TDA into machine learning workflows to improve classification accuracy and interpretability. Cawi et al. (2019) explored the incorporation of Mapper into predictive modeling, demonstrating its effectiveness in structured feature extraction. Their study applied TDA Mapper approaches to hospital readmission and credit evaluation tasks, where traditional data pipelines often lacked interpretability and flexibility. The results showed that integrating Mapper into machine learning workflows significantly improved classification accuracy. Building on these advancements, Riihimäki et al. (2020) developed a TDA Mapper classification method for repeated measurements, addressing a common challenge in machine learning models that struggle with multiple measurement points per sample. Their method outperformed standard SVM models, achieving 96.8% accuracy in point process classification and 90% in tree species classification, while also offering insights into high-purity data subgroups, making it particularly valuable for biological sciences.

To improve the resolution parameter selection in the Mapper algorithm, Kang and Lim (2021) proposed an ensemble-based approach to enhance the stability and reliability of Mapper outputs. Their study focused on two resolution parameters, the number of intervals and the overlapping percentage, which significantly impact the consistency of data representations. To address this issue, they introduced an ensemble Mapper method that integrates multiple Mapper results generated under different parameter settings using a fuzzy clustering ensemble technique, ensuring more robust and consistent topological representations. Their findings demonstrated that the ensemble method was more robust and stable compared to conventional single-parameter selection methods, significantly improving visualization and clustering accuracy in high-dimensional data.

Most recently, De Lara (2023) further refined Mapper's role in classification by integrating it into persistent homology-based classification frameworks. The study demonstrated that Mapper-enhanced classification models outperformed traditional approaches like Support Vector Machines (SVM) and Random Forest, offering improved scalability and predictive performance. These collective advancements illustrate the growing role of TDA in machine learning, showcasing its ability to enhance classification accuracy, optimize parameter selection, and improve the interpretability of complex models. Samani & Banerjee (2024) extended TDA's applications to robotics, introducing THOR2, a topological analysis framework for 3D shape and color-based human-inspired object recognition. Conventional deep learning models often struggle with domain generalization in unseen environments. By integrating the Mapper algorithm with color-based embeddings, their model outperformed deep learning approaches on real-world datasets, improving robotic perception in cluttered indoor settings.

### *3.1.5 Biomedicine and Genomics*

Yao et al. (2009) were among the earliest to apply TDA in the biological sciences, specifically to biomolecular folding pathways. Conventional clustering methods often obscure transient states in protein folding due to their rarity in simulation data. By applying the Mapper algorithm, Yao et al. successfully identified multiple intermediate states in RNA hairpin folding, providing new insights into folding mechanisms. A few years later, Lum et al. (2013) explored the application of TDA for extracting geometric patterns in complex datasets. Traditional statistical methods like PCA often failed to detect subtle yet biologically significant subgroups. Their study demonstrated how TDA could uncover previously undetected clusters in breast cancer gene expression, U.S. congressional voting behavior, and NBA player performance statistics, revealing new insights for targeted therapies, policy analysis, and sports analytics.

In the following years, TDA continued to gain traction in biomedical and genomic research. Campbell et al. (2019) applied the TDA Mapper for feature extraction in pain recognition using physiological signals, addressing the challenges of inconsistent feature selection in emotion recognition systems. Around the same time, Wang et al. (2019) tackled challenges in single-cell RNA sequencing (scRNA-seq) visualization. Traditional methods like t-SNE struggled to preserve continuous gene expression trajectories, whereas Mapper effectively captured clustering structures through topological graphs and maintained gene expression topologies, aiding in lineage differentiation studies. Feng et al. (2019) further contributed by introducing a semi-supervised topological analysis method (STA) for transcriptomic data. Their findings demonstrated that Mapper-generated topological graphs were more effective than PCA in distinguishing biologically relevant structures. Additionally, Amézquita et al. (2020) applied TDA to developmental biology, uncovering novel morphological transformations across species through the Mapper output topological graph.

As research expanded into clinical applications, Lopez et al. (2022) applied TDA using Mapper Plus to classify patients based on their risk of adverse events (AEs) following acute myocardial infarction (AMI). Analyzing 798 AMI patients through 31 clinical variables, TDA identified three distinct subgroups with varying survival probabilities, providing a hypothesis-free, data-driven approach to clinical risk stratification. This method offered a new framework for automated risk assessment without requiring prior clinical assumptions. In the same year, Amézquita et al. (2023) integrated the TDA Mapper with spectral shape analysis to analyze RNA sequencing data from lung adenocarcinoma samples. Traditional clustering methods often overlooked subtle biological differences in tumors, but this approach successfully identified two distinct subgroups of tumor cells with unique gene expression patterns, revealing potential pathways for cancer progression through the topological graph.

More recent studies have continued refining the role in medical and genomic research. Ma et al. (2023) applied the Mapper algorithm to analyze the progression of motor features in Parkinson's disease, constructing a Markov chain on Mapper topological graphs to model motor deterioration over time. Their results allowed patient-specific predictions of motor decline, providing valuable insights for clinicians adjusting treatment strategies. In plant biology, Palande et al. (2023) explored how gene expression shapes plant morphology using TDA Mapper. Their study identified a core expression backbone defining form and function across diverse plant lineages, offering a novel computational framework for studying plant evolution. Percival et al. (2024) extended these applications by investigating the development of leaf shape in grapevine and maracuyá. Standard statistical methods failed to detect conserved heteroblastic patterns across species, whereas TDA Mapper, through the topological graph, effectively captured evolutionary and developmental signatures hidden in high-dimensional morphological data.

Finally, Mapper has demonstrated its impact in predictive modeling and diagnostics. A recent study in 2024 applied TDA with the Mapper algorithm to predict heart disease using two UCI heart disease datasets (Cleveland and Statlog). Traditional machine learning models often struggled with feature selection and dimensionality reduction, resulting in suboptimal performance. To overcome these challenges, researchers applied a tri-dimensional Singular Value Decomposition (SVD) filter, improving predictive accuracy to an impressive 99.32% on the Cleveland dataset and 99.62% on the Statlog dataset. Loughrey et al. (2024) highlight the Mapper algorithm's role in patient subgroup discovery by capturing topological and geometric structures in biomedical data. However, traditional machine learning struggles with classifying new patients in small Mapper graph regions. To address this, they introduce the Mapper k-Nearest Neighbor (Mapper-kNN) algorithm, which improves classification using graph-based distances. Experiments on five gene expression cancer datasets confirm its effectiveness for more accurate patient classification. These results reinforce the importance of the TDA Mapper in medical diagnostics and healthcare analytics, further solidifying its role as a powerful tool in high-dimensional data analysis.

### *3.1.6 Physics and Planetary Science*

Murugan and Robertson (2019) introduced the TDA Mapper to the field of physics, emphasizing its advantages over traditional clustering methods. They applied persistent homology and the Mapper algorithm to fast radio burst (FRB) data, demonstrating that the Mapper topological graph could uncover meaningful patterns in astrophysical observations that conventional statistical techniques often failed to detect. Their study showcased how topology-based methods can contribute to the understanding of complex, high-dimensional physical data. Lukmanov et al. (2021) applied the TDA Mapper to laser ablation ionization mass spectrometry (fs-LIMS) data for planetary exploration. Identifying biosignatures in planetary materials is challenging due to the complexity of spectral data. Their analysis of Precambrian Gunflint microfossils showed that the Mapper algorithm could effectively classify distinct mineralogical and organic signatures through topological graphs, demonstrating its potential for future space missions.

### *3.1.7 Epidemiology and Public Health*

Chen & Volić (2021) utilized the TDA Mapper algorithm to model the spread of COVID-19 in the U.S., addressing a key limitation of conventional epidemiological models, which struggle to integrate spatial and temporal dynamics simultaneously. By constructing Mapper topological graphs of pandemic data, the study provided a more comprehensive visualization of COVID-19 transmission patterns, effectively revealing growth rates and regional hotspots that traditional statistical models often overlooked. Around the same time, Kitanishi et al. (2021) applied TDA to analyze depression patient journeys using claims data. Traditional epidemiological models often fail to capture disease progression across multiple clinical outcomes. The study leveraged TDA Mapper topological graph visualization to classify diseases associated with depression, enhancing precision medicine approaches for pharmaceutical interventions.

TDA Mapper's role in improving disease modeling, medical diagnostics, and patient stratification continued to expand. Rafique and Mir (2020) developed a TDA Mapper approach for cancer subtyping using gene expression data. Their Mapper-based framework significantly enhanced subtype differentiation and survival separability, identifying relevant genetic pathways critical for precision oncology. Dagliati et al. (2020) further applied the Mapper algorithm to disease progression modeling in type 2 diabetes. By integrating pseudo-time series analysis with Mapper, they identified distinct disease trajectories from electronic health records (EHRs), highlighting lipid profiles as key factors in individualized treatment predictions. In cardiovascular research, Hwang et al. (2021) classified coronary atherosclerosis progression using coronary computed tomography angiography (CTA). Their study identified three distinct patient groups with varying plaque dynamics, showing that individuals with moderate plaque burden had the highest incidence of acute coronary syndrome. Mapper-derived patient classifications significantly improved cardiovascular risk prediction beyond conventional methods.

More recent studies have continued to refine TDA applications in disease modeling. Rather and Chachoo (2023) tackled the challenge of identifying clinically relevant disease subtypes using omics data, which is inherently sparse and high-dimensional. Their robust pipeline, integrating UMAP with Mapper, outperformed traditional clustering techniques by achieving superior separation in survival outcomes across multiple cancer datasets. Chao et al. (2023) applied the Mapper to identify prognostically distinct phenotypes in patients undergoing transcatheter edge-to-edge repair (TEER). Their model, which incorporated mitral valve morphology and clinical risk factors, revealed two major patient clusters with significantly different survival outcomes, allowing for better risk stratification in TEER procedures.

Ling et al. (2024) used the Mapper algorithm to analyze COVID-19 data in Malaysia from January 2020 to June 2022, identifying major pandemic waves and regional variations. It highlights Selangor as an outlier due to high case numbers and shows distinct patterns in Sabah and Sarawak. Using clustering techniques and feature selection, the study compares Mapper with hierarchical clustering and PCA, demonstrating its effectiveness in capturing complex epidemiological trends. While sensitive to parameter selection, the approach provides a compact and effective tool for visualizing pandemic progression and informing public health decisions. These advancements demonstrate the growing impact of the TDA Mapper in biomedical research, offering novel methods to analyze complex disease patterns, enhance precision medicine, and improve patient outcome predictions.

### *3.1.8 Cardiovascular Analysis*

Tokodi et al. (2020) explored patient similarity networks in cardiac function to enhance personalized cardiovascular medicine. Traditional analytical methods often struggle to integrate multiple echocardiographic features for predicting major adverse cardiac events (MACE). To address this challenge, the study applied unsupervised machine learning and TDA Mapper to develop a phenotypic network or topological graph for classifying patients based on cardiac function. Their model improved risk stratification and provided a framework for the automated staging of cardiac disease severity, offering a more comprehensive approach to cardiovascular diagnostics.

Building on these advancements, Hernández-Lemus et al. (2024) provided an extensive review of TDA applications in cardiovascular signal analysis, focusing on the limitations of traditional methods in handling high-dimensional and heterogeneous data. By leveraging persistent homology and the Mapper algorithm, the study demonstrated how TDA can significantly enhance cardiovascular disease (CVD) diagnosis and prognosis. Their findings showed that TDA-based approaches allow for a more detailed analysis of electrocardiography and arterial stiffness, leading to improved patient outcomes and more precise clinical decision-making. These studies underscore the growing role of TDA in cardiovascular medicine, highlighting its potential to refine diagnostics and improve risk assessment strategies

### *3.1.9 Economics and Finance*

Guo and Banerjee (2017) were among the first to apply the Mapper algorithm in manufacturing system predictions, demonstrating its effectiveness in process monitoring. Their study identified key variables in chemical process yield prediction and semiconductor wafer fault detection while maintaining high predictive accuracy. By leveraging Mapper, they improved process efficiency and defect detection, providing a more robust framework for manufacturing system analysis. In the field of economics and technological innovation, Escolar et al. (2019) applied the Mapper algorithm to patent data to visualize firms' locations in technological space. Their study addressed the challenges of mapping high-dimensional and unstructured innovation trajectories. The results revealed that firms with unique patenting behaviors, represented as "flares" in the Mapper topological graph, were linked to higher financial performance, offering new insights into how firms innovate and compete over time. Qiu et al. (2020) extended the application of the Mapper into corporate failure prediction by integrating it with Altman’s Z-score model. Their findings showed that failing firms were not always clustered in distinct regions of financial ratio space, challenging traditional assumptions about bankruptcy risk assessment and providing a new perspective on corporate financial health.

As research into Mapper’s economic applications progressed, Escolar et al. (2023) expanded their earlier work by exploring firms’ innovation trajectories through patent data, further refining the visualization of technological space. Their study demonstrated that firms with unique inventive trajectories tended to perform better, uncovering previously undocumented industry dynamics. In finance, Shiraj et al. (2024) tackled anomaly detection in financial time series data. Standard clustering techniques like DBSCAN often failed to detect subtle financial anomalies. By combining Mapper with DBSCAN clustering, the study provided a more refined anomaly detection framework, improving financial risk prediction and market stability assessments. Uray et al. (2024) review topological data analysis in smart manufacturing, categorizing applications into quality control, process monitoring, and manufacturing engineering. Using persistent homology, Mapper, and UMAP, they analyze industrial data for anomaly detection, material flow optimization, and predictive maintenance. While TDA shows promise, its adoption remains limited due to challenges in interpretability and scalability.

Looking ahead, Wei et al. (2025) examined the role of Mapper and topology in economic and financial modeling. Traditional econometric models struggle with high-dimensional financial data and complex market structures. Their study demonstrated that Mapper and topological machine learning techniques improve risk assessment, asset correlation analysis, and financial time series prediction, highlighting the growing influence of TDA in modern financial analytics. Collectively, these studies underscore the increasing relevance of Mapper across manufacturing, innovation analysis, corporate risk assessment, and financial modeling, showcasing its versatility in tackling complex data challenges.

### 3.1.10 Education

Anderson et al. (2022) applied topological data analysis (TDA) to analyze the Serbia 2014–2015 Multiple Indicator Cluster Surveys (MICS) dataset, which provides household data on wealth, health, and education. Using the Mapper algorithm, they transformed the data into a topological graph, uncovering relationships between household wealth and factors such as urban/rural living, ownership of assets, and prioritization of possessions. Their findings suggest that TDA can provide insights into socio-economic disparities and inform policy decisions, potentially refining the wealth index and highlighting key determinants of household well-being. Restrepo (2023) applies topological data analysis to track grade evolution in university math courses before, during, and after COVID-19. Using 24 datasets on grades, subjects, and teacher details, the study employs two TDA algorithms to analyze fail rate scenarios. Findings show that the Mapper algorithm effectively identifies significant grade differences based on teachers' contract types.

### 3.1.11 Cheminformatics and Material Sciences

TDA has also been applied in cheminformatics to enhance molecular modeling and prediction accuracy. Pirashvili et al. (2018) tackled the limitations of traditional solubility prediction models, which often suffered from overfitting and lacked interpretability. Using the Mapper algorithm, the study analyzed chemical descriptor space and uncovered unexpected relationships between molecular features such as chlorine content and solubility. The findings demonstrated that topological methods could improve solubility prediction by identifying novel descriptors and enhancing model interpretability.

### 3.1.12 Biological Sciences

TDA Mapper has played a transformative role in medical and genomic research by uncovering patterns in high-dimensional biological data. Mathews et al. (2019) applied the Mapper algorithm to the molecular phenotyping of soft tissue sarcomas, addressing the complexity of genomic datasets. By combining network geometry, diffusion maps, and the Mapper algorithm, the study effectively identified distinct gene expression signatures corresponding to various sarcoma subtypes, offering deeper insights into their molecular heterogeneity. Notably, a novel signature linked to TP53 tumor suppressor activation was discovered, offering a new method for classifying molecular subtypes in cancer research. Jeitziner et al. (2019) further contributed to gene expression analysis by developing Two-Tier Mapper (TTMap), a topology-based clustering method that improved the stability and sensitivity of clustering techniques in high-dimensional biological datasets, allowing for more reliable identification of significant subgroups.

As the Mapper applications in medicine expanded, Carr et al. (2021) introduced a topological machine learning approach to identify homogeneous patient subgroups in medical data. Traditional clustering methods struggled with mixed data types and non-linear relationships, but by leveraging the Mapper algorithm, the study successfully grouped patients with shared prognostic outcomes, advancing precision medicine. Their pipeline also incorporated prior knowledge, enhancing robustness and interpretability. More recently, Imoto and Hiraoka (2023) introduced V-Mapper, an extension of the Mapper algorithm designed to integrate velocity information in single-cell RNA sequencing (scRNA-seq) data. This approach addressed challenges in modeling dynamic cellular processes and differentiation pathways, successfully mapping biologically meaningful cell-state transitions by embedding RNA velocity into the Mapper topological graph.

Beyond medicine, Mapper has also influenced the study of biological shapes and structures. Amézquita et al. (2020) applied TDA to analyze morphological structures at molecular and organismal levels, demonstrating its ability to capture and quantify shape variations that traditional statistical methods often overlook. Their findings provided a novel framework for studying biological shape changes over time. More recently, Derwae et al. (2023) applied Mapper to mass spectrometry imaging data, successfully identifying spatial molecular distributions in pancreatic tissue. Their study not only replicated known structures but also discovered novel biological patterns, such as blood vessel clustering, showcasing Mapper's ability to uncover hidden biological relationships. These studies highlight the growing influence of TDA across diverse scientific fields. From enhancing clustering stability in gene expression data and advancing neuroimaging analysis to improving disease modeling and biological classification, TDA Mapper continues to be a powerful tool for extracting meaningful insights from complex, high-dimensional datasets.

## 3.2 TDA Mapper Algorithm Variants

Over the years, several variants of the Mapper algorithm have been developed to enhance its efficiency, stability, and versatility, making it more adaptable to a wide range of complex data analysis tasks. These advancements have improved clustering accuracy, parameter selection, computational efficiency, and visualization capabilities, enabling more robust and interpretable representations of high-dimensional datasets. The continuous evolution of these Mapper variants has expanded its applicability across diverse domains, including biology, medicine, finance, artificial intelligence, and engineering, further solidifying its role as a powerful tool for uncovering meaningful patterns in complex data, as discussed in the following section.

### *3.2.1 Ball Mapper*

Recently, Dłotko (2019a) introduced Ball Mapper, an algorithm inspired by the conventional Mapper algorithm, designed to enhance the topological representation of high-dimensional data. While Ball Mapper shares the goal of revealing the shape of complex datasets, it is generally considered a distinct approach rather than a direct variant of Mapper due to its different construction principles. Unlike conventional Mapper, which relies on multiple parameters such as cover size, resolution, and clustering to construct a topological graph, Ball Mapper simplifies this process by employing a collection of closed balls of equal radius to cover all data points within a dataset. Each ball serves as a node containing a subset of similar data points, and edges are formed between nodes that share common members. This construction results in a Ball Mapper topological graph, structurally similar to conventional Mapper outputs, but with a more intuitive and parameter-reduced formation method. The nodes in Ball Mapper function as clusters of similar data points, making it particularly useful for understanding local structures in complex datasets.

A key advantage of Ball Mapper over conventional Mapper is its reduced parameter dependency, as it requires only a single parameter, the radius of the ball, whereas conventional Mapper demands multiple parameters to generate the topological graph. This makes Ball Mapper easier to implement, though it might be more difficult to interpret, especially in cases where selecting appropriate parameters in a conventional Mapper can be challenging. While both methods produce rich topological graphs capable of capturing structural relationships within a dataset, there has been limited research comparing their effectiveness across various applications. Despite the relatively scarce literature on Ball Mapper, existing studies demonstrate its potential in revealing topological structures and enabling qualitative data analysis. Current applications of Ball Mapper include COVID-19 clinical data analysis (Dłotko & Rudkin, 2020) and financial analysis (Dłotko et al., 2021; Dłotko et al., 2022), illustrating its ability to extract meaningful insights from diverse datasets.

More recently, Valerio et al. (2024) applied Ball Mapper to study the three-dimensional radiodensity distribution of fish otoliths. Understanding otolith density variations is crucial for studying fish metabolism, growth, and population dynamics. Their study optimized sampling methods for dimensionality reduction while maintaining biologically relevant density variations, validating the use of Ball Mapper for reducing computational costs in high-dimensional biological data analysis. Similarly, Rudkin et al. (2024) demonstrate that the predictability of Bitcoin returns depends significantly on past return trajectories, with specific patterns yielding more accurate forecasts using Ball Mapper. Madukpe et al. (2025) conducted a comparative study of Ball Mapper and conventional Mappers to assess air pollutant behavior across 60 monitoring stations in Malaysia. Their findings demonstrated that while both methods provided valuable insights into air quality trends, Ball Mapper's simplicity in requiring only a single parameter made it more effective for analyzing pollutant patterns and detecting redundancies in monitoring stations. This research suggests Ball Mapper as a preferable tool for air quality monitoring and regulatory policy improvements.

Finally, Park et al. (2025) utilized Ball Mapper to examine spatial interactions between malignant and immune cells within the tumor microenvironment, leveraging high-dimensional imaging data from diffuse large B-cell lymphoma (DLBCL). By constructing topological graphs based on 18-dimensional features encoding distance and density metrics, Ball Mapper offered a powerful framework for visualizing complex spatial patterns and immune infiltration. This approach enabled the identification of subtle transitions between tumor-enriched and border regions, addressing limitations inherent to traditional local-centric analyses. Significantly, Ball Mapper graphs were combined with hierarchical clustering and survival analysis, uncovering spatial phenotypes linked to patient outcomes. This demonstrated its potential as a valuable tool for identifying prognostic biomarkers in heterogeneous tumor landscapes.

### 3.2.2 F-Mapper

The Fuzzy-Mapper (F-Mapper) algorithm, introduced by Bui et al. (2020), is an extension of the Mapper algorithm designed to improve the automatic division of cover intervals in TDA. Unlike the original Mapper, which requires predefined interval length and overlap percentage parameters, F-Mapper utilizes Fuzzy c-Means (FCM) clustering to dynamically determine these intervals. This approach eliminates rigid interval divisions and allows for more natural, data-driven coverings, making it more adaptable to complex and high-dimensional datasets. The study tested F-Mapper on three datasets, Unit Circle, Reaven and Miller Diabetes, and NKI Breast Cancer demonstrated that it produces topologically similar structures to Mapper while improving clustering quality, as measured by silhouette coefficient scores.

The F-Mapper algorithm retains the core advantages of Mapper, such as its ability to extract meaningful patterns from high-dimensional data, but introduces greater flexibility and automation in parameter selection. By employing Fuzzy c-Means clustering, F-Mapper ensures that cover intervals vary based on data density rather than following a uniform structure. This results in a more robust, adaptive, and interpretable topological representation, particularly in biomedical and machine learning applications. The study concludes that F-Mapper enhances clustering accuracy, reduces computational constraints associated with traditional Mapper implementations, and provides a more effective tool for qualitative data analysis. Future research is expected to further refine its applications in various fields, particularly in bioinformatics, neuroscience, and social sciences.

### 3.2.3 Ensemble-based Mapper

Kang and Lim (2021) present an enhanced approach to the Mapper algorithm by addressing its sensitivity to two resolution parameters, specifically the number of intervals and the overlap percentage. They propose an ensemble technique that generates multiple Mapper outputs under varying parameter settings and integrates them using fuzzy clustering. This method demonstrates improved robustness and performance across three real-world datasets from the fields of biology and chemistry. Similarly, Fitzpatrick et al. (2023) is the first work to address the challenge of selecting three core parameters of the Mapper algorithm simultaneously, namely the lens function, the resolution, and the gain. Their approach employs ensemble learning to construct multiple Mapper graphs across a wide range of parameter values and identifies the most stable and representative structures. The method involves three main steps: first, selecting optimal cover parameters for each lens function; second, identifying common data structures across different lenses; and third, constructing an ensemble Mapper graph by combining the most representative graphs through sample co-occurrence clustering.

### 3.2.4 Two-Tier Mapper (TTMap)

Jeitziner et al. (2019) introduced Two-Tier Mapper (TTMap), a variant designed to improve clustering stability in high-dimensional biological data. TTMap integrates hierarchical and partitioning clustering techniques, making it more robust and less dependent on user-defined parameters compared to traditional clustering methods.

### 3.2.5 V-Mapper

Imoto and Hiraoka (2023) developed Velocity Mapper (V-Mapper), a variant of the Mapper algorithm that incorporates velocity information in single-cell RNA sequencing (scRNA-seq) data. V-Mapper uses RNA velocity as the filter function within the Mapper framework to capture the dynamic trajectories of cells over time. By embedding RNA velocity into the Mapper graph construction process, it enables the accurate mapping of cell-state transitions, offering deeper insights into cellular differentiation and dynamic biological processes. Once the data is segmented according to the RNA velocity filter values, standard clustering methods are applied within each segment to group similar cells, allowing the final Mapper topological graph to represent the structural relationships and temporal progression inherent in the data.

### 3.2.6 G-Mapper

Alvarado et al. (2023) introduced G-Mapper, an enhanced variant of the traditional Mapper algorithm that focuses on optimizing the cover construction. Unlike the traditional approach, which relies on a predefined number of intervals with uniform length and overlap, G-Mapper employs the G-means clustering algorithm to adaptively and automatically partition the data. This allows the cover to better reflect the underlying data distribution, resulting in more meaningful and stable Mapper graphs. It starts with a single cover interval and iteratively splits it by testing whether the data inside follows a Gaussian distribution using the Anderson-Darling test. If the data doesn't fit a normal distribution, the interval is split into two, using a Gaussian Mixture Model (GMM) to guide where and how to split. This approach allows G-Mapper to create a data-driven, non-uniform cover that better reflects the true structure of the dataset, improving the interpretability and stability of the resulting Mapper graph. Unlike traditional Mapper methods, G-Mapper does not require the user to predefine the number of intervals or their overlap; instead, it learns these parameters automatically from the data

#### 3.2.7 *Mapper-Enhanced Persistent Homology Classification Algorithm (MEPHCA)*

De Lara and Garcia (2024) developed the Mapper-Enhanced Persistent Homology Classification Algorithm (MEPHCA), a variant that focuses on improving computational efficiency in classification tasks. MEPHCA reduces the traditionally high computational costs associated with Mapper while maintaining classification accuracy, making it more practical for large-scale applications.

#### 3.2.8 *D-Mapper*

Tao and Ge (2025) introduced a distribution-guided Mapper algorithm (D-Mapper), an enhanced version of the conventional Mapper algorithm used in TDA. Unlike the traditional approach, which relies on fixed interval lengths and overlapping ratios that can miss subtle data features, D-Mapper leverages probabilistic models- specifically, mixture models to generate data-driven, flexible covers that better reflect the data's intrinsic distribution. The authors propose a novel evaluation metric that integrates clustering quality, assessed using the Silhouette Coefficient, with topological structure analysis, measured through Extended Persistent Homology. This combined approach provides a more comprehensive assessment of Mapper algorithm outputs, ensuring both well-defined cluster separation and the preservation of meaningful topological features. Experimental results, including applications to synthetic datasets and SARS-CoV-2 RNA sequences, demonstrate that D-Mapper consistently outperforms the classic Mapper, offering deeper insights into both vertical and horizontal evolutionary processes.

### 3.3 Strengths and Potential Weaknesses of the Mapper Algorithm Variants

The various Mapper algorithm variants represent efforts to address key limitations inherent in the conventional Mapper framework. Traditional Mapper requires the manual selection of multiple sensitive parameters lens function, cover intervals, resolution, gain, and clustering method, making it prone to inconsistent outputs, subjective tuning, and unstable interpretations, especially when analyzing complex, high-dimensional datasets. Ball Mapper simplifies this by reducing parameter dependence to a single radius value, streamlining implementation while retaining the ability to capture local topological features. However, Ball Mapper's simplicity can come at the cost of interpretability, as choosing an appropriate radius lacks a formal guideline and can significantly affect the resulting graph. On the other hand, F-Mapper introduces a data-driven approach by leveraging fuzzy clustering to dynamically adjust covers, improving adaptability and reducing reliance on predefined parameters. While this approach enhances adaptability by reducing reliance on predefined parameters, it also removes the interpretability provided by a filter function. In conventional Mapper, the filter function allows users to explore different topological perspectives of the data, whereas in F-Mapper, the clustering dynamics alone dictate the final structure. This can make it more difficult to control or interpret the features being captured in the resulting topological graph. The ensemble-based Mapper addresses the challenges of parameter selection in conventional Mappers by generating and aggregating multiple Mapper graphs across varying filter functions, covering parameters and clustering methods. It employs graph merging techniques and stability measures to synthesize these graphs into a unified representation. However, while this approach enhances robustness, the extensive computation time required for an exhaustive search of parameter combinations to construct instability matrices increases computational complexity. TTMap improves clustering stability by combining hierarchical and partitioning methods. Still, its two-level structure could make it harder to scale to larger datasets and more difficult to interpret the results.

Further developments such as V-Mapper, G-Mapper, D-Mapper, and MEPHCA extend Mapper's scope into specific analytical challenges and data modalities. V-Mapper incorporates temporal dynamics through RNA velocity as a filter function, allowing it to capture time-dependent changes in data. However, its accuracy depends heavily on the reliability of the RNA velocity estimates, which can be affected by data quality and noise. G-Mapper addressed parameter sensitivity with statistical rigor, yet it depends on assumptions of normality that may not hold in all datasets. D-Mapper enhances cover generation through probabilistic modeling, offering better sensitivity to data distribution, although probabilistic approaches can add model complexity and require careful validation. MEPHCA focuses on computational efficiency in classification but might trade off general-purpose flexibility. Across these innovations, the persistent gap in conventional Mappers lies in balancing parameter automation, computational efficiency, and interpretability; no single variant fully resolves all three simultaneously. Each approach mitigates specific limitations of the conventional Mapper, yet often introduces trade-offs in complexity, scalability, or application scope.

## 3.4 Limitations of the Mapper Algorithm and Future Research Directions

Despite the various variants of the Mapper algorithm developed to enhance the conventional Mapper, several limitations remain unaddressed. The variants only focused on the choice of cover and clustering method specific to a field, while Ball Mapper is the only variant that is completely independent of filter function and clustering algorithm. These gaps present opportunities for future research to explore and improve upon. Developing robust statistical frameworks for Mapper analysis and automated parameter optimization techniques.

1. *Sensitivity to Parameter Selection*: One of the major challenges of the Mapper algorithm is its high sensitivity to parameter choices, including the filter function, covering parameters, and clustering method. Even slight variations in these parameters can result in significantly different Mapper graphs, making the outcomes unstable and difficult to reproduce. This sensitivity introduces a degree of subjectivity, as different parameter choices can emphasize different structures within the data, leading to inconsistent interpretations and reducing the reliability of the analysis.

*Future Research Directions:* Adaptive parameter selection, including machine learning and optimization techniques, could be developed to dynamically adjust Mapper parameters based on the characteristics of the data, although Fitzpatrick et al. (2023) have introduced an automatic parameter selection technique, yet more are needed. Automated hyperparameter tuning methods, such as Bayesian optimization, genetic algorithms, or reinforcement learning, could help identify the most stable and meaningful parameter configurations. Finally, robust sensitivity analysis could be improved by developing quantitative measures to assess the stability and reliability of the Mapper graph across different parameter settings.

2. *Dependence on the Filter Function*: The Mapper algorithm relies heavily on the choice of the filter function, which determines how the dataset is projected into a lower-dimensional space before applying covering and clustering. This dependence presents several challenges: (a) Different filter functions can produce vastly different Mapper graphs, leading to inconsistencies in the results. (b) There is no universally accepted method for selecting the optimal filter function, making the choice subjective and data-dependent. (c) Some filter functions may distort the topology of the data, resulting in the loss of meaningful structures such as clusters, loops, or voids.

*Future Research Directions*: Machine learning models could be developed to automatically determine the most suitable filter function based on the characteristics of the dataset. Additionally, using multiple filter functions simultaneously and combining their results could lead to more robust and consistent representations. Finally, designing new filter functions that preserve intrinsic topological features, such as loops and voids, would help minimize structural distortions and improve the accuracy of Mapper-based analyses.

3. *Potential Loss of Global Structure*: One fundamental limitation of the Mapper algorithm is its tendency to simplify high-dimensional data into a network representation, which may lead to the loss of critical global topological and geometric structures. This issue arises primarily due to two factors: (a) Mapper relies on a lens function to reduce the dataset's dimensionality and partitions it into local regions before clustering, which can disrupt global connectivity patterns, and (b) if the lens function and partitioning strategy are not carefully designed, essential structures such as cycles, holes, or continuous manifolds may be lost.

*Future Research Directions*: To address these challenges, future research should focus on enhancing Mapper's ability to retain global topological features. One promising approach is integrating the Mapper algorithm with persistent homology, which could enable the tracking and preservation of critical structures like loops and connected components across different scales. Persistent homology ensures that key topological properties remain intact, mitigating distortions caused by clustering and dimensionality reduction. Additionally, developing data-dependent covering methods can improve the accuracy of Mapper graphs by ensuring they retain meaningful structural properties, ultimately leading to more robust and interpretable representations of complex datasets.

4. *Interpretability Challenges and Limited Theoretical Guarantees*: While Mapper provides a graph-based summary of data, interpreting its output is not always straightforward. Several factors contribute to this challenge: (a) different parameter settings, such as lens functions, cover partitions, and clustering methods, can yield significantly different Mapper graphs, making it difficult to draw definitive conclusions; (b) the lack of standardized quantitative evaluation metrics makes it challenging to assess the quality and reliability of Mapper representations; and (c) there are no generally acceptable rigorous theoretical bounds to formally relate Mapper outputs to the true underlying topology of the dataset, although some progress has been made in this regard by Carrière et al. (2017) raising concerns about consistency and robustness.

*Future Research Directions*: To enhance the interpretability and theoretical foundations of Mapper, future research should focus on developing quantitative measures that assess the quality of Mapper outputs based on topological stability and consistency; however, there is recent progress in this regard by Tao and Ge (2025), but more is needed. Additionally, establishing more rigorous theoretical frameworks that guarantee the robustness and reliability of Mapper graphs will be essential for wider adoption in scientific and engineering applications.

# 4 Conclusion

The present review comprehensively reports the advances in the development, refinement, and application of the Mapper algorithm, a key technique within topological data analysis (TDA). Through an analysis of peer-reviewed literature spanning from 2007 to 2025, this study highlights how the Mapper algorithm has been instrumental in visualizing and interpreting complex, high-dimensional datasets across diverse fields, such as neuroscience, biomedicine, genomics, environmental science, finance, and machine learning. It provides a valuable resource for researchers and practitioners, offering insights into methodological innovations, including emerging Mapper variants like Ball Mapper, F-Mapper, and V-Mapper, which address key challenges related to scalability, stability, and parameter sensitivity. From the comparison of the lens functions on two concentric circle-shaped point cloud datasets, PCA and MDS perform better, while GMM, TDA-GC, HACA, and DBSCAN perform well as clustering algorithms within the Mapper algorithm by producing topological graphs that best represent the structure of the dataset.

It is demonstrated in this review that the Mapper algorithm holds immense potential as an advanced data analysis tool capable of revealing meaningful insights from complex datasets. Its ability to simplify high-dimensional data into intuitive topological graphs makes it an indispensable technique for modern data science, particularly in the era of big data production. The emergence of user-friendly software implementations and automated frameworks has further expanded its accessibility and applicability. Despite existing challenges, the Mapper algorithm continues to play a pivotal role in uncovering novel insights and supporting informed decision-making in a wide range of scientific, industrial, and technological fields. Its adoption and versatility are expected to enhance the quality and depth of data-driven discoveries and applications in the years to come.

# Acknowledgments

The authors deeply thank Universiti Sains Malaysia for their financial support through a Short-Term Grant (Project No: 304/PMATHS/6315649).

# Statements and Declarations

## Funding

Universiti Sains Malaysia supported this work through a Short-Term Grant (Project No: 304/PMATHS/6315649).

### Competing Interests

The authors have no competing interests to declare relevant to this article's content.

### Author Contributions

Conceptualization: Vine Nwabuisi Madukpe; Methodology: Vine Nwabuisi Madukpe and Bright Chukwuma Ugoala; Formal analysis and investigation: Vine Nwabuisi Madukpe and Bright Chukwuma Ugoala; Fariha Syaqina Zulkepli; Writing - original draft preparation: Vine Nwabuisi Madukpe and Bright Chukwuma Ugoala. Writing - review and editing: Vine Nwabuisi Madukpe and Bright Chukwuma Ugoala. Supervision: Fariha Syaqina Zulkepli..

### Ethical Approval

This declaration is not applicable, as this study does not involve human images, human data, or the use of animals.

# References

Anderson, J. R., Memić, F., & Volić, I. (2022). Topological data analysis and UNICEF multiple indicator cluster surveys. *Journal of Quantitative Economics*, *20*(2), 281–309. https://doi.org/10.1007/s40953-022-00288-w

Amézquita, E. J., Nasrin, F., Storey, K. M., & Yoshizawa, M. (2023). Genomics data analysis via spectral shape and topology. *PLoS ONE*, *18*(4), e0284820. https://doi.org/10.1371/journal.pone.0284820

Amézquita, E. J., Quigley, M. Y., Ophelders, T., Munch, E., & Chitwood, D. H. (2020). The shape of things to come: Topological data analysis and biology, from molecules to organisms. *Developmental Dynamics*, *249*(7), 816–833. https://doi.org/10.1002/dvdy.175

Almgren, K., Kim, M., & Lee, J. (2017). Extracting Knowledge from the Geometric Shape of Social Network Data Using Topological Data Analysis. *Entropy*, *19*(7), 360. https://doi.org/10.3390/e19070360

Aljanobi, F. A., & Lee, J. (2021). Topological data analysis for classification of heart disease data. *In Proceedings of the 2021 IEEE International Conference on Big Data and Smart Computing (BigComp)* (pp. 210–213). IEEE.

Ayasdi. (n.d.). *Ayasdi AI platform*. Retrieved February 27, 2025, from https://www.ayasdi.com/

Bui, Q., Vo, B., Nguyen, H., DO, Hung, N. Q. V., & Snasel, V. (2019). F-Mapper: a fuzzy Mapper clustering algorithm. *Knowledge-Based Systems*, *189*, 105107. https://doi.org/10.1016/j.knosys.2019.105107

Belchi, F., Brodzki, J., Burfitt, M., & Niranjan, M. (2020). A numerical measure of the instability of Mapper-Type algorithms. *arXiv (Cornell University)*, *21*(202), 1–45. https://arxiv.org/pdf/1906.01507.pdf

Bungula, W., & Darcy, I. (2024). Bi-Filtration and stability of TDA Mapper for point cloud data. *arXiv (Cornell University)*. https://doi.org/10.48550/arxiv.2409.17360

Chen, Y., & Volić, I. (2021). Topological data analysis model for the spread of the coronavirus. *PLoS ONE*, *16*(8), e0255584. https://doi.org/10.1371/journal.pone.0255584

Carr, E., Carrière, M., Michel, B., Chazal, F., & Iniesta, R. (2021). Identifying homogeneous subgroups of patients and important features: a topological machine learning approach. *BMC Bioinformatics*, *22*(1). https://doi.org/10.1186/s12859-021-04360-9

Côté-Allard, U., Campbell, E., Phinyomark, A., Laviolette, F., Gosselin, B., & Scheme, E. (2020). Interpreting Deep learning features for Myoelectric Control: A comparison with Handcrafted features. *Frontiers in Bioengineering and Biotechnology*, *8*. https://doi.org/10.3389/fbioe.2020.00158

Campbell, E., Phinyomark, A., & Scheme, E. (2019). Feature extraction and selection for pain recognition using peripheral physiological signals. *Frontiers in Neuroscience*, *13*. https://doi.org/10.3389/fnins.2019.00437

Chao, C., Barry, T., Seri, A., Shaer, A. E., Ponce, N. C., Chakraborty, S., Smith, S., Alkhouli, M., Thaden, J., Fortuin, D., Sweeney, J. P., Eleid, M., Rihal, C. S., Holmes, D. R., Pollak, P. M., Sabbagh, A. E., Lester, S. J., Oh, J. K., Shen, W., . . . Arsanjani, R. (2023). Topological data analysis identified Prognostically-Distinct phenotypes in transcatheter Edge-to-Edge repair patients. *Mayo Clinic Proceedings Digital Health*, *1*(3), 381–392. https://doi.org/10.1016/j.mcpdig.2023.07.002

Cawi, E., La Rosa, P. S., & Nehorai, A. (2019). Designing machine learning workflows with an application to topological data analysis. *PLoS ONE*, *14*(12), e0225577. https://doi.org/10.1371/journal.pone.0225577

Carrière, M., Michel, B., & Oudot, S. (2017). Statistical analysis and parameter selection for Mapper. *arXiv (Cornell University)*. https://doi.org/10.48550/arxiv.1706.00204

Chazal, F., & Michel, B. (2021). An Introduction to Topological Data Analysis: Fundamental and Practical Aspects for Data Scientists. *Frontiers in Artificial Intelligence*, *4*. https://doi.org/10.3389/frai.2021.667963

Carlsson, G. (2020). Topological methods for data modelling. *Nature Reviews Physics*, *2*(12), 697–708. https://doi.org/10.1038/s42254-020-00249-3

Corcoran, P., & Jones, C. B. (2023). Topological data analysis for geographical information science using persistent homology. *International Journal of Geographical Information Science*, *37*(3), 712–745. https://doi.org/10.1080/13658816.2022.2155654

Darling, R. W. R., Emanuello, J. A., Purvine, E., & Ridley, A. (2022). Proceedings of TDA: Applications of Topological Data Analysis to Data Science, Artificial Intelligence, and Machine Learning Workshop at SDM 2022. *arXiv (Cornell University)*. https://doi.org/10.48550/arxiv.2204.01142

Daou, L., & Hanna, E. M. (2024). Predicting protein complexes in protein interaction networks using Mapper and Graph Convolution Networks. *Computational and Structural Biotechnology Journal*, *23*, 3595–3609. https://doi.org/10.1016/j.csbj.2024.10.009

Derwae, H., Nijs, M., Geysels, A., Waelkens, E., & De Moor, B. (2023). Spatiochemical characterization of the pancreas using mass spectrometry imaging and topological data analysis. *Analytical Chemistry*, *95*(28), 10550–10556. https://doi.org/10.1021/acs.analchem.2c05606

Duman, A. N., & Tatar, A. E. (2023). Topological data analysis for revealing dynamic brain reconfiguration in MEG data. *PeerJ*, *11*, e15721. https://doi.org/10.7717/peerj.15721

Duman, A. N., Tatar, A. E., & Pirim, H. (2019). Uncovering dynamic brain reconfiguration in MEG working Memory N-Back task using topological data analysis. *Brain Sciences*, *9*(6), 144. https://doi.org/10.3390/brainsci9060144

Dłotko, P., Gurnari, D., & Sazdanovic, R. (2024). Mapper–Type algorithms for complex data and relations. *Journal of Computational and Graphical Statistics*, *33*(4), 1383–1396. https://doi.org/10.1080/10618600.2024.2343321

Dłotko, P., Qiu, W., & Rudkin, S. T. (2021). Financial ratios and stock returns reappraised through a topological data analysis lens. *The European Journal of Finance*, 1–25. https://doi.org/10.1080/1351847x.2021.2009892

Dłotko, P. (2019a). Ball mapper: a shape summary for topological data analysis. *ArXiv (Cornell University)*. https://doi.org/10.48550/arxiv.1901.07410

Dlotko, P. (2019b). *BallMapper: The Ball Mapper algorithm* [Dataset]. https://doi.org/10.32614/cran.package.ballmapper

Dlotko, P. (2019c) The Ball Mapper Algorithm *[R package BallMapper version 0.2.0]*. https://cran.r-project.org/package=BallMapper

Dlotko, P., & Rudkin, S. (2020). Covid-19 clinical data analysis using Ball Mapper. *MedRxiv (Cold Spring Harbor Laboratory)*. https://doi.org/10.1101/2020.04.10.20061374

Dłotko, P., Qiu, W., & Rudkin, S. (2022). Topological Data Analysis Ball Mapper for Finance. *ArXiv (Cornell University)*. https://doi.org/10.48550/arxiv.2206.03622

Dagliati, A., Geifman, N., Peek, N., Holmes, J. H., Sacchi, L., Bellazzi, R., Sajjadi, S. E., & Tucker, A. (2020). Using topological data analysis and pseudo time series to infer temporal phenotypes from electronic health records. *Artificial Intelligence in Medicine*, *108*, 101930. https://doi.org/10.1016/j.artmed.2020.101930

Escolar, E. G., Hiraoka, Y., Igami, M., & Ozcan, Y. (2023b). Mapping firms' locations in technological space: A topological analysis of patent statistics. *Research Policy*, *52*(8), 104821. https://doi.org/10.1016/j.respol.2023.104821

Endo, S., & Yokogawa, S. (2021). Analysis of the trends between indoor carbon dioxide concentration and Plug-Level electricity usage through topological data analysis. *IEEE Sensors Journal*, *22*(2), 1424–1434. https://doi.org/10.1109/jsen.2021.3130570

Fan, L., Li, Y., Huang, Z., Zhang, W., Wu, X., Liu, T., & Wang, J. (2023). Low-frequency repetitive transcranial magnetic stimulation alters the individual functional dynamical landscape. *Cerebral Cortex*, *33*(16), 9583–9598. https://doi.org/10.1093/cercor/bhad228

Faaique, M. (2023). Overview of Big data analytics in modern astronomy. *International Journal of Mathematics Statistics and Computer Science*, *2*, 96–113. https://doi.org/10.59543/ijmscs.v2i.8561

Fitzpatrick, P., Jurek-Loughrey, A., Dlotko, P., & Rincon, J. M. D. (2023). *Ensemble learning for Mapper parameter optimization. In 2023 IEEE 35th International Conference on Tools with Artificial Intelligence (ICTAI) (pp. 129–134). IEEE.* https://doi.org/10.1109/ICTAI59109.2023.00026

Feng, T., Davila, J. I., Liu, Y., Lin, S., Huang, S., & Wang, C. (2019). Semi-Supervised Topological Analysis for elucidating hidden structures in High-Dimensional Transcriptome datasets. *IEEE/ACM Transactions on Computational Biology and Bioinformatics*, *18*(4), 1620–1631. https://doi.org/10.1109/tcbb.2019.2950657

Foughrey, C. F., Maguire, S., Dłotko, P., Bai, L., Orr, N., & Jurek-Loughrey, A. (2025). A novel method for subgroup discovery in precision medicine based on topological data analysis. *BMC Medical Informatics and Decision Making*, 25, 139. https://doi.org/10.1186/s12911-025-02852-9

Gurnari, P. (2022). *pyBallMapper: Python implementation of the Ball Mapper algorithm* [Software]. GitHub. https://github.com/Gurnari/pyBallMapper

Guo, W., & Banerjee, A. G. (2017a). Identification of key features using topological data analysis for accurate prediction of manufacturing system outputs. *Journal of Manufacturing Systems*, *43*, 225–234. https://doi.org/10.1016/j.jmsy.2017.02.015

Geniesse, C., Sporns, O., Petri, G., & Saggar, M. (2019). Generating dynamical neuroimaging spatiotemporal representations (DyNeuSR) using topological data analysis. *Network Neuroscience*, *3*(3), 763–778. https://doi.org/10.1162/netn_a_00093

Geniesse, C., Chowdhury, S., & Saggar, M. (2022). NeuMapper: A scalable computational framework for multiscale exploration of the brain's dynamical organization. *Network Neuroscience*, 1–32. https://doi.org/10.1162/netn_a_00229

Haşegan, D., Geniesse, C., Chowdhury, S., & Saggar, M. (2024). Deconstructing the Mapper algorithm to extract richer topological and temporal features from functional neuroimaging data. *Network Neuroscience*, *8*(4), 1355–1382. https://doi.org/10.1162/netn_a_00403

Hwang, D., Kim, H. J., Lee, S., Lim, S., Koo, B., Kim, Y., Kook, W., Andreini, D., Al-Mallah, M. H., Budoff, M. J., Cademartiri, F., Chinnaiyan, K., Choi, J. H., Conte, E., Marques, H., De Araújo Gonçalves, P., Gottlieb, I., Hadamitzky, M., Leipsic, J. A., . . . Chang, H. (2021). Topological data analysis of coronary plaques demonstrates the natural history of coronary atherosclerosis. *JACC. Cardiovascular Imaging*, *14*(7), 1410–1421. https://doi.org/10.1016/j.jcmg.2020.11.009

Hernández-Lemus, E., Miramontes, P., & Martínez-García, M. (2024). Topological Data Analysis in Cardiovascular Signals: An Overview. *Entropy*, *26*(1), 67. https://doi.org/10.3390/e26010067

Imoto, Y., & Hiraoka, Y. (2023). V-Mapper: topological data analysis for high-dimensional data with velocity. *Nonlinear Theory and Its Applications IEICE*, *14*(2), 92–105. https://doi.org/10.1587/nolta.14.92

Islambekov, U., Yuvaraj, M., & Gel, Y. R. (2019b). Harnessing the power of topological data analysis to detect change points. *Environmetrics*, *31*(1). https://doi.org/10.1002/env.2612

Jeitziner, R., Carrière, M., Rougemont, J., Oudot, S., Hess, K., & Brisken, C. (2019). Two-Tier Mapper, an unbiased topology-based clustering method for enhanced global gene expression analysis. *Bioinformatics*, *35*(18), 3339–3347. https://doi.org/10.1093/bioinformatics/btz052

Kaveh-Yazdy, F., & Zarifzadeh, S. (2020). Track Iran's national COVID-19 response committee's major concerns using two-stage unsupervised topic modeling. *International Journal of Medical Informatics*, *145*, 104309. https://doi.org/10.1016/j.ijmedinf.2020.104309

Kyeong, S., Park, S., Cheon, K., Kim, J., Song, D., & Kim, E. (2015). A New Approach to Investigate the Association between Brain Functional Connectivity and Disease Characteristics of Attention-Deficit/Hyperactivity Disorder: Topological Neuroimaging Data Analysis. *PLoS ONE*, *10*(9), e0137296. https://doi.org/10.1371/journal.pone.0137296

Kitanishi, Y., Fujiwara, M., & Binkowitz, B. (2021). Patient journey through cases of depression from claims database using machine learning algorithms. *PLoS ONE*, *16*(2), e0247059. https://doi.org/10.1371/journal.pone.0247059

Kang, S. J., & Lim, Y. (2021). Ensemble mapper. *Stat*, *10*(1). https://doi.org/10.1002/sta4.405

Kim, D., Kim, S., Kim, Y. D., & Lyu, S. (2025). Enhanced detection of harmful algal blooms using topological data analysis for clustering spatially distributed water quality and hydrodynamic data. *KSCE Journal of Civil Engineering*, 100177. https://doi.org/10.1016/j.kscej.2025.100177

Lauric, A., Ludwig, C. G., & Malek, A. M. (2023). Topological data analysis and use of Mapper for cerebral aneurysm rupture status discrimination based on 3-Dimensional shape analysis. *Neurosurgery*, *93*(6), 1285–1295. https://doi.org/10.1227/neu.0000000000002570

Ling, C. Y., Phang, P., Liew, S., Jayaraj, V. J., & Wiwatanapataphee, B. (2024b). Exploration of COVID-19 data in Malaysia through mapper graph. *Network Modeling Analysis in Health Informatics and Bioinformatics*, *13*(1). https://doi.org/10.1007/s13721-024-00472-3

Liao, T., Wei, Y., Luo, M., Zhao, G., & Zhou, H. (2019). tmap: an integrative framework based on topological data analysis for population-scale microbiome stratification and association studies. *Genome Biology*, *20*(1). https://doi.org/10.1186/s13059-019-1871-4

Lum, P. Y., Singh, G., Lehman, A., Ishkanov, T., Vejdemo-Johansson, M., Alagappan, M., Carlsson, J., & Carlsson, G. (2013b). Extracting insights from the shape of complex data using topology. *Scientific Reports*, *3*(1). https://doi.org/10.1038/srep01236

Loughrey, C. F., Dłotko, P., & Jurek-Loughrey, A. (2024). A Mapper-Based classifier for patient subgroup prediction. In *IFMBE proceedings* (pp. 610–621). https://doi.org/10.1007/978-3-031-62502-2_69

Lopez, J. E., Datta, E., Ballal, A., & Izu, L. T. (2022). Abstract 14875: Topological Data Analysis of Electronic Health Record Features Predicts major cardiovascular outcomes after revascularization for acute myocardial infarction. *Circulation*, *146*(Suppl_1). https://doi.org/10.1161/circ.146.suppl_1.14875

Lukmanov, R. A., Riedo, A., Wacey, D., Ligterink, N. F. W., Grimaudo, V., Tulej, M., De Koning, C., Neubeck, A., & Wurz, P. (2021). On Topological analysis of FS-LIMS data. Implications for in situ planetary mass spectrometry. *Frontiers in Artificial Intelligence*, *4*. https://doi.org/10.3389/frai.2021.668163

Lucasimi. (n.d.). *tda-mapper GitHub repository*. Retrieved February 27, 2025, from https://github.com/lucasimi/tda-mapper-python

Madukpe, V. N., Zulkepli, N. F. S., Noorani, M. S. M., & Gobithaasan, R. U. (2025). Comparative analysis of Ball Mapper and conventional Mapper in investigating air pollutants' behavior. *Environmental Monitoring and Assessment*, *197*(2). https://doi.org/10.1007/s10661-024-13477-2

Munch, E. (2017b). A user's guide to topological data analysis. *Journal of Learning Analytics*, *4*(2). https://doi.org/10.18608/jla.2017.42.6

Ma, L., Feng, T., He, C., Li, M., Ren, K., & Tu, J. (2023). A progression analysis of motor features in Parkinson's disease based on the mapper algorithm. *Frontiers in Aging Neuroscience*, *15*. https://doi.org/10.3389/fnagi.2023.1047017

Ma, G. (2020). Using topological data analysis to process time-series data: a persistent homology way. *Journal of Physics Conference Series*, *1550*(3), 032082. https://doi.org/10.1088/1742-6596/1550/3/032082

Mathews, J. C., Pouryahya, M., Moosmüller, C., Kevrekidis, I., Deasy, J. O., & Tannenbaum, A. (2018). Molecular phenotyping using networks, diffusion, and topology: soft tissue sarcoma. *bioRxiv (Cold Spring Harbor Laboratory)*. https://doi.org/10.1101/328054

Müllner, D. and Babu, A. (2013). *Python Mapper: An open-source toolchain for data exploration, analysis, and visualization.* URL http://danifold.net/mapper

Musa, S. M. S., Md Noorani, M. S., Abdul Razak, F., Ismail, M., Alias, M. A., & Hussain, S. I. (2021). Using persistent homology as preprocessing of early warning signals for critical transition in

flood. *Scientific Reports*, *11*(1). https://doi.org/10.1038/s41598-021-86739-5
Nielson, J. L., Paquette, J., Liu, A. W., Guandique, C. F., Tovar, C. A., Inoue, T., Irvine, K., Gensel, J. C., Kloke, J., Petrossian, T. C., Lum, P. Y., Carlsson, G. E., Manley, G. T., Young, W., Beattie, M. S., Bresnahan, J. C., & Ferguson, A. R. (2015). Topological data analysis for discovery in preclinical spinal cord injury and traumatic brain injury. *Nature Communications*, *6*(1). https://doi.org/10.1038/ncomms9581
Nicolau, M., Levine, A. J., & Carlsson, G. (2011). Topology based data analysis identifies a subgroup of breast cancers with a unique mutational profile and excellent survival. *Proceedings of the National Academy of Sciences*, *108*(17), 7265–7270. https://doi.org/10.1073/pnas.1102826108
Palande, S., Kaste, J. a. M., Roberts, M. D., Abá, K. S., Claucherty, C., Dacon, J., Doko, R., Jayakody, T. B., Jeffery, H. R., Kelly, N., Manousidaki, A., Parks, H. M., Roggenkamp, E. M., Schumacher, A. M., Yang, J., Percival, S., Pardo, J., Husbands, A. Y., Krishnan, A., . . . VanBuren, R. (2023). Topological data analysis reveals a core gene expression backbone that defines form and function across flowering plants. *PLoS Biology*, *21*(12), e3002397. https://doi.org/10.1371/journal.pbio.3002397
Percival, S., Onyenedum, J. G., Chitwood, D. H., & Husbands, A. Y. (2024). Topological data analysis reveals core heteroblastic and ontogenetic programs embedded in leaves of grapevine (Vitaceae) and maracuyá (Passifloraceae). *PLoS Computational Biology*, *20*(2), e1011845. https://doi.org/10.1371/journal.pcbi.1011845
Pirashvili, M., Steinberg, L., Guillamon, F. B., Niranjan, M., Frey, J. G., & Brodzki, J. (2018). Improved understanding of aqueous solubility modeling through topological data analysis. *Journal of Cheminformatics*, *10*(1). https://doi.org/10.1186/s13321-018-0308-5
Paultpearson. (2015). GitHub - paultpearson/TDAmapper: (R package) Analyze High-Dimensional Data Using Discrete Morse Theory. *GitHub*. https://github.com/paultpearson/TDAmapper/
Patania, A., Selvaggi, P., Veronese, M., Dipasquale, O., Expert, P., & Petri, G. (2019). Topological gene expression networks recapitulate brain anatomy and function. *Network Neuroscience*, *3*(3), 744–762. https://doi.org/10.1162/netn_a_00094
Pearson, P. T. (2013). Visualizing clusters in artificial neural networks using MORSE theory. *Advances in Artificial Neural Systems*, *2013*, 1–8. https://doi.org/10.1155/2013/486363
Restrepo, M., & Cornelis, C. (2024). Mapper-Based rough sets. In *Lecture notes in computer science* (pp. 3–17). https://doi.org/10.1007/978-3-031-65665-1_1
Restrepo, M. (2023). Topological data analysis for the evolution of student grades before, during and after the COVID-19 pandemic. In *Studies in big data* (pp. 97–119). https://doi.org/10.1007/978-3-031-38325-0_5
Pedoia, V., Haefeli, J., Morioka, K., Teng, H., Nardo, L., Souza, R. B., Ferguson, A. R., & Majumdar, S. (2017). MRI and biomechanics multidimensional data analysis reveals R2-R1ρ as an early predictor of cartilage lesion progression in knee osteoarthritis. *Journal of Magnetic Resonance Imaging*, *47*(1), 78–90. https://doi.org/10.1002/jmri.25750
Qiu, W., Rudkin, S., & Dłotko, P. (2020). Refining understanding of corporate failure through a topological data analysis mapping of Altman's Z-score model. *Expert Systems With Applications*, *156*, 113475. https://doi.org/10.1016/j.eswa.2020.113475
Rather, A. A., & Chachoo, M. A. (2023). Robust correlation estimation and UMAP assisted topological analysis of omics data for disease subtyping. *Computers in Biology and Medicine*, *155*, 106640. https://doi.org/10.1016/j.compbiomed.2023.106640
Rafique, O., & Mir, A. (2019). A topological approach for cancer subtyping from gene expression data. *Journal of Biomedical Informatics*, *102*, 103357. https://doi.org/10.1016/j.jbi.2019.103357
Robles, A., Hajij, M., & Rosen, P. (2018). The shape of an image - a study of Mapper on Images. *Proceedings of the 17th International Joint Conference on Computer Vision, Imaging and Computer Graphics Theory and Applications*, 339–347. https://doi.org/10.5220/0006574803390347
Riihimäki, H., Chachólski, W., Theorell, J., Hillert, J., & Ramanujam, R. (2020). A topological data analysis based classification method for multiple measurements. *BMC Bioinformatics*, *21*(1). https://doi.org/10.1186/s12859-020-03659-3
Rudkin, S., Rudkin, W., & Dłotko, P. (2024). Return trajectory and the forecastability of bitcoin returns. *Financial Review*. https://doi.org/10.1111/fire.12420
Szairis. (n.d.). *GitHub - szairis/sakmapper: Implementation of Mapper Algorithm*. GitHub. https://github.com/szairis/sakmapper
Simi, L. (2024). tda-mapper (v0.8.0). Zenodo. https://doi.org/10.5281/zenodo.14194667
Singh, G., Mémoli, F., & Carlsson, G. E. (2007). Topological methods for the analysis of high dimensional data sets and 3D object recognition. *Eurographics*, 91–100. https://doi.org/10.2312/spbg/spbg07/091-100
Singh, Y., Hathaway, Q. A., Farrelly, C., Budoff, M. J., Erickson, B., Collins, J. D., Blaha, M. J., Leiner, T., Lopez-Jimenez, F., Rozenblit, J., Sarkar, D., & Carlsson, G. (2025). Topological Data analysis in the Assessment of Coronary Atherosclerosis: A Comprehensive Narrative review. *Mayo Clinic Proceedings Digital Health*, 100199. https://doi.org/10.1016/j.mcpdig.2025.100199

Saggar, M., Sporns, O., Gonzalez-Castillo, J., Bandettini, P. A., Carlsson, G., Glover, G., & Reiss, A. L. (2018). Towards a new approach to reveal dynamical organization of the brain using topological data analysis. *Nature Communications*, *9*(1). https://doi.org/10.1038/s41467-018-03664-4

Saggar, M., Bruno, J., Gaillard, C., Claudino, L., & Ernst, M. (2022a). Neural resources shift under Methylphenidate: A computational approach to examine anxiety-cognition interplay. *NeuroImage*, *264*, 119686. https://doi.org/10.1016/j.neuroimage.2022.119686

Saggar, M., Shine, J. M., Liégeois, R., Dosenbach, N. U. F., & Fair, D. (2022b). Precision dynamical mapping using topological data analysis reveals a hub-like transition state at rest. *Nature Communications*, *13*(1). https://doi.org/10.1038/s41467-022-32381-2

Singer, B., Meling, D., Hirsch-Hoffmann, M., Michels, L., Kometer, M., Smigielski, L., Dornbierer, D., Seifritz, E., Vollenweider, F. X., & Scheidegger, M. (2024). Psilocybin enhances insightfulness in meditation: a perspective on the global topology of brain imaging during meditation. *Scientific Reports*, *14*(1). https://doi.org/10.1038/s41598-024-55726-x

Samani, E. U., & Banerjee, A. G. (2024). THOR2: Topological analysis for 3D shape and Color-Based Human-Inspired object recognition in unseen environments. *Advanced Intelligent Systems*. https://doi.org/10.1002/aisy.202400539

Sun, C. (2020). Exploration of Mapper-A method for Topological data Analysis. *2020 International Conference on Information Science, Parallel and Distributed Systems (ISPDS)*, 142–145. https://doi.org/10.1109/ispds51347.2020.00036

Savic, A., Toth, G., & Duponchel, L. (2017). Topological data analysis (TDA) applied to reveal pedogenetic principles of European topsoil system. *The Science of the Total Environment*, *586*, 1091–1100. https://doi.org/10.1016/j.scitotenv.2017.02.095

Skaf, Y., & Laubenbacher, R. (2022b). Topological data analysis in biomedicine: A review. *Journal of Biomedical Informatics*, *130*, 104082. https://doi.org/10.1016/j.jbi.2022.104082

Shah, W. H., Jaimes-Reátegui, R., Huerta-Cuellar, G., García-López, J., & Pisarchik, A. (2025). Persistent homology approach for uncovering transitions to Chaos. *Chaos Solitons & Fractals*, *192*, 116054. https://doi.org/10.1016/j.chaos.2025.116054

Siddiqui, S., Shikotra, A., Richardson, M., Doran, E., Choy, D., Bell, A., Austin, C. D., Eastham-Anderson, J., Hargadon, B., Arron, J. R., Wardlaw, A., Brightling, C. E., Heaney, L. G., & Bradding, P. (2018). Airway pathological heterogeneity in asthma: Visualization of disease microclusters using topological data analysis. *Journal of Allergy and Clinical Immunology*, *142*(5), 1457–1468. https://doi.org/10.1016/j.jaci.2017.12.982

Tokodi, M., Shrestha, S., Bianco, C., Kagiyama, N., Casaclang-Verzosa, G., Narula, J., & Sengupta, P. P. (2020). Interpatient similarities in cardiac function. *JACC. Cardiovascular Imaging*, *13*(5), 1119–1132. https://doi.org/10.1016/j.jcmg.2019.12.018

Todd, J. T., & Petrov, A. A. (2022). The many facets of shape. *Journal of Vision*, *22*(1), 1. https://doi.org/10.1167/jov.22.1.1

Tao, Y., Ge, S. A. (2025). Distribution-guided Mapper algorithm. *BMC Bioinformatics* 26, 73. https://doi.org/10.1186/s12859-025-06085-5

van Veen, H., Saul, N., Eargle, D., & Mangham, S. (2019). Kepler Mapper: A flexible Python implementation of the Mapper algorithm. *Journal of Open Source Software*, *4*(42), 1315. https://doi.org/10.21105/joss.01315

Vannoni, S., Tassi, E., Sampaio, I. W., Bianchi, A. M., & Maggioni, E. (2024). A Systematic approach to tuning cover parameters in Mapper for improved TDA representation. *2022 IEEE International Conference on Metrology for Extended Reality, Artificial Intelligence and Neural Engineering (MetroXRAINE)*, 377–382. https://doi.org/10.1109/metroxraine62247.2024.10796460

Valerio, J., Vasconcelos-Filho, J. E., Stosic, B., De Oliveira, W. R., Santana, F. M., Antonino, A. C. D., & Duarte-Neto, P. J. (2024). Topological Analysis of the Three-Dimensional radiodensity distribution of fish otoliths: point sampling effects on dimensionality Reduction. *Micron*, *188*, 103731. https://doi.org/10.1016/j.micron.2024.103731

Wang, T., Johnson, T., Zhang, J., & Huang, K. (2019). Topological methods for visualization and analysis of high dimensional Single-Cell RNA sequencing data. *PMC*. https://scholarworks.iupui.edu/bitstream/1805/20958/1/nihms-999822.pdf

Walsh, K., Voineagu, M. A., Vafaee, F., & Voineagu, I. (2020). TDAview: an online visualization tool for topological data analysis. *Bioinformatics*, *36*(18), 4805–4809. https://doi.org/10.1093/bioinformatics/btaa600

Yao, Y., Sun, J., Huang, X., Bowman, G. R., Singh, G., Lesnick, M., Guibas, L. J., Pande, V. S., & Carlsson, G. (2009). Topological methods for exploring low-density states in biomolecular folding pathways. *The Journal of Chemical Physics*, *130*(14). https://doi.org/10.1063/1.3103496

Zhang, M., Chowdhury, S., & Saggar, M. (2022). Temporal Mapper: Transition networks in simulated and real neural dynamics. *Network Neuroscience*, *7*(2), 431–460. https://doi.org/10.1162/netn_a_00301

Zhou, Y., Jenne, H., Brown, D., Shapiro, M., Jefferson, B., Joslyn, C., Henselman-Petrusek, G., Praggastis, B., Purvine, E., & Wang, B. (2023). Comparing mapper graphs of artificial neuron activations. *Proceedings - 2023 Topological Data Analysis and Visualization, TopoInVis 2023*, *34*, 41–50. https://doi.org/10.1109/topoinvis60193.2023.00011

Zhou, Y., Chalapathi, N., Rathore, A., Zhao, Y., & Wang, B. (2020). Mapper Interactive: a scalable, extendable, and interactive toolbox for the visual exploration of High-Dimensional data. *arXiv (Cornell University)*. https://doi.org/10.48550/arxiv.2011.03209

Zulkepli, N. F. S., Noorani, M. S. M., Razak, F. A., Ismail, M., & Alias, M. A. (2022). Hybridization of hierarchical clustering with persistent homology in assessing haze episodes between air quality monitoring stations. *Journal of Environmental Management*, *306*, 114434. https://doi.org/10.1016/j.jenvman.2022.114434

Zulkepli, N. F. S., Madukpe, V. N., Noorani, M. S. M., Bakar, M. a. A., Gobithaasan, R. U., & Jie, O. C. (2024). Topological clustering in investigating spatial patterns of particulate matter between air quality monitoring stations in malaysia. *Air Quality Atmosphere & Health*. https://doi.org/10.1007/s11869-024-01596-1